\documentclass[aap,MSNbibl,dvips]{arximspdf}
\usepackage{ifthen}
\usepackage{graphicx}

\doi{10.1214/12-AAP871} %kopijuoti is PTS
\volume{23}
\issue{4}
\pubyear{2013}
\firstpage{1318}
\lastpage{1354}

\makeatletter
\newcommand{\rrVert}{\Vert}
\newcommand{\rrvert}{\vert}
\newcommand{\llVert}{\Vert}
\newcommand{\llvert}{\vert}

\newcommand{\lldots}{,\ldots,}

\newcommand{\s}[1]{a_{#1}}
\newcommand{\equilib}{\mathcal E}
\newcommand{\eq}{\equilib}
\newcommand{\bd}[1]{\dot{\bar{#1}}}
\newcommand{\f}[1]{\ifthenelse{\equal{#1}{l}}{\mathrm{f}_{(l)}}{\mathrm{f}_{(#1)}}}

\newcommand{\ff}[1]{\f{#1}}
\newcommand{\R}{\mathbb{R}}
\newcommand{\E}{\mathrm{E}}

\newcommand{\Ni}{\mathcal{N}}
\newcommand{\Pb}{\mathrm{P}}

\newtheorem{theorem}{Theorem}
\newtheorem{lemma}{Lemma}

\newproclaim{definition}{Definition}

\newtheorem{corollary}{Corollary}

\newproclaim{property}{Property}
\newproclaim{fact}{Fact}

\makeatother

\begin{document}
\begin{frontmatter}

\title{Fluid limits to analyze long-term flow rates
of a~stochastic network with ingress discarding\thanksref{T1}}
\runtitle{Flow rate analysis using fluid limits}

\thankstext{T1}{Supported in part by NSF Grants ANI-0331659, CNS-0953884
and CNS-0910702.}

\begin{aug}
\author[A]{\fnms{John} \snm{Musacchio}\corref{}\ead[label=e1]{johnm@soe.ucsc.com}}
\and
\author[B]{\fnms{Jean} \snm{Walrand}\ead[label=e2]{wlr@eecs.berkeley.edu}}
\runauthor{J. Musacchio and J. Walrand}
\affiliation{University of California, Santa Cruz, and University of
California, Berkeley}
\address[A]{Technology and Information Management\\
University of California, Santa Cruz \\
1156 High Street\\
MS: SOE 3\\
Santa Cruz, California 95064 \\
USA\\
\printead{e1}}
\address[B]{Department of Electrical Engineering\\
\quad and Computer Sciences \\
University of California, Berkeley\\
257M Cory Hall \\
Berkeley, California 94720 \\
USA\\
\printead{e2}} %adresu isvedimo komanda gale!
\end{aug}

% HISTORY:
\received{\smonth{3} \syear{2010}}
\revised{\smonth{4} \syear{2012}}

% ABSTRACT
%
\begin{abstract}
We study a simple rate control scheme for a multiclass queuing network
for which customers are partitioned into distinct flows that are queued
separately at each station. The control scheme discards customers that
arrive to the network ingress whenever any one of the flow's queues
throughout the network holds more than a specified threshold number of
customers. We prove that if the state of a corresponding fluid model
tends to a set where the flow rates are equal to target rates, then
there exist sufficiently high thresholds that make the long-term
average flow rates of the stochastic network arbitrarily close to these
target rates. The same techniques could be used to study other control
schemes. To illustrate the application of our results, we analyze a
network resembling a 2-input, 2-output communications network switch.
% to show that the long-term average flow rates approach those of the
%fluid model.
\end{abstract}

% KEYWORDS

\end{frontmatter}

%s1 #&#
\section{Introduction}
\label{sintro}

We consider a multiclass queuing network whose customers are
partitioned into $F$ distinct flows. Customers of a flow $f
\in
\{1,\ldots,F\}$ arrive according to an independent renewal process and
follow a fixed, acyclic sequence of stations. The service times at each
station are also independent. Each flow $f$ has a weight $w_f\in\R_+$,
and each of $d$ stations is equipped with per-flow queues and serves a
flow in proportion to its weight using a weighted round robin or a
similar queueing discipline like weighted fair queueing or generalized
head of line processor sharing.

We consider a simple scheme which we call ingress discarding for
admitting customers. The ingress discarding scheme works as follows.
Whenever any of a flow's queues exceed a threshold $h$, that flow's
customers are discarded at the network ingress. There are two main
objectives of the scheme: (i) stability when the arrival rates in the
absence of discarding would cause the utilization of some stations to
exceed 1, and (ii) fairness in the long-term average departure rates
when the network cannot accommodate all the incoming flows. The
contribution of this article is a methodology for proving that the
long-term average flow rates in such a network can be made arbitrarily
close to those predicted by a fluid model, provided that the discarding
thresholds are sufficiently high.

\label{applications} There are a number of applications of such a
control policy. One application is for service centers such as call
centers. It might be acceptable to block incoming customers, but
unacceptable to drop customers that have been admitted to the system,
hence the appropriateness of ingress discarding. A designer of such a
system might want to show that the flow rates of various types of
customers are fair in some sense. This work can be used to show that if
the system's fluid model achieves fair rates, then the system will
achieve close to fair rates provided that the discarding thresholds are
sufficiently high. Another application area is in data-packet switch
design. A packet switch typically consists of several line-cards that
transmit and receive the data packets, and a switch-fabric that serves
as an interconnect. A design requirement might be that any packet
discarding occur in the line-cards rather than in the switch fabric,
since the line cards are better equipped to record statistics about the
dropped packets, for instance. The switch fabric can be thought of as a
queuing network, and ingress discarding would be one way to fulfill the
requirement that discarding only occur in the line cards. Again, this
work shows that the flow rates of such a system approach those
predicted by a fluid model if the discarding thresholds are made
sufficiently high.

\begin{figure}

\includegraphics{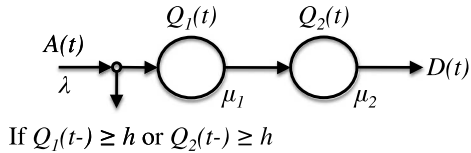}

\caption{A queueing network with input discarding.}
\label{figex1}
\end{figure}

To illustrate our methodology, we consider the simple network in Figure
\ref{figex1}.
%
%f1 #&#
%
This network carries a single flow and customers arrive as a renewal
process $E(t)$. There are two queues, each with i.i.d. service times
with mean $\mu_i^{-1}$ in queue $i$ ($i = 1, 2$). Designate by $Q_i(t)$
the length of queue $i$ ($i = 1, 2$). The ingress discarding scheme
discards the arrivals that occur when one of the two lengths is at
least equal to threshold $h$. We want to show that if the thresholds
are made large enough factor~$n$, that the flow rates approach $\min\{
\lambda, \mu_1, \mu_2 \}$. More precisely, we want to show that for
every $\epsilon> 0$ there exists some $n_\epsilon$ such that if
threshold scale factor $n \geq n_\epsilon$, then the average rate of
the departure process $D(t)$ exceeds $\min\{\lambda, \mu_1, \mu_2
\}
- \epsilon$. \label{Note}Note that since we scale the thresholds by a
factor $n$, the starting value of the threshold $h$ is not important,
so long as it is positive. Also note that we do not attempt to derive
any result on the speed of convergence---how fast $nh$ must grow to
achieve rates within a smaller and smaller $\epsilon$ of the desired rates.

\begin{figure}

\includegraphics{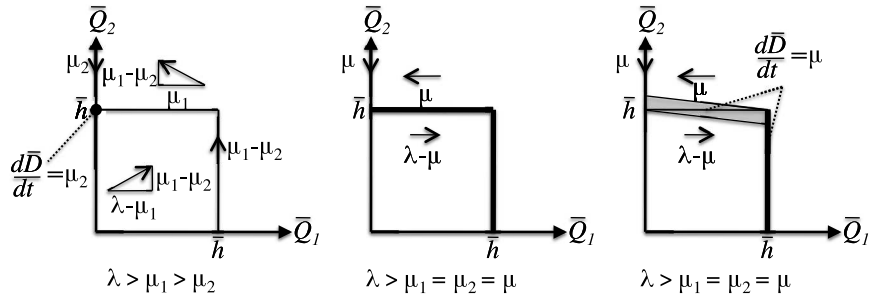}

\caption{The fluid process that approximates the stochastic network.}
\label{figex2}
\end{figure}

\label{approach} The analysis approach, which we believe can be
extended to control strategies that change admission, service, or
routing behavior when queue depths cross thresholds that can be made
large, is based on deriving properties of the stochastic network using
a fluid model. However for clarity of exposition, we limit our focus in
this paper to the ingress discarding policy. As in work by Dai \cite
{Dai95} we take a fluid limit by considering a sequence of larger and
larger initial conditions, and scaling time and space by the size of
those initial conditions. However, in order to consider stochastic
networks with larger and larger thresholds, our fluid limit also
considers a sequence of systems with thresholds scaled by an increasing
factor $n$. The resulting fluid limit behaves according to a fluid
model corresponding to the vector flow diagram in Figure~\ref{figex2}.
Since we scale the thresholds in our fluid limit, the thresholds appear
in the fluid model with nonnegligible values $\bar{h}$. Note that
$\bar
{h}$ need not equal $h$ since the fluid limits we consider may scale
space and threshold at different rates. Also as a consequence of
scaling the thresholds in taking the fluid limit, the stochastic system
behaves like the fluid model (in terms of flow rates) only if the
stochastic system's thresholds are sufficiently large.

%
%f2 #&#
%

First consider the case $\lambda> \mu_1 > \mu_2$. A fluid model
corresponding to this case is illustrated by the vector flow diagram in
the left part of Figure~\ref{figex2}.
This diagram indicates the rate of change of the vector of queue
lengths as a function of its value. For instance, if the two queue
lengths are between 0 and $\bar{h}$, then fluid enters queue 1 at rate
$\lambda$ and flows from that queue to queue 2 at rate $\mu_1$ while
fluid leaves queue 2 at rate $\mu_2$. Accordingly, the length of queue
1 increases at rate $\lambda- \mu_1$ and that of queue 2 at rate $\mu_1
- \mu_2$. The other cases can be understood similarly. The vector
flow diagram shows that, irrespective of their initial values, the
queue lengths converge to the pair of values $(0, \bar{h})$, which is
an absorbing state for the fluid process. Moreover, when the process is
close to the value $(0, \bar{h})$, the rate of the departure fluid is
close to $\mu_2$.
To conclude that the stochastic network has a departure rate close to
$\mu_2$ when $h$ is large, one notes that the fluid process has one
additional property: the time the process takes to reach the state $(0,
\bar{h})$ is bounded by a linear function of the distance between the
initial condition and $(0, \bar{h})$. This property, which can be seen
from the vector flow diagram, can be used to show, roughly, that the
stochastic system spends little time far from $(0, \bar{h})$. The
intuition is that, although fluctuations occasionally move the
stochastic network away from the limiting state, the system tends to
follow the fluid process and get back to that state fairly quickly.
This property will allow us to construct a proof that the stochastic
network has a departure rate close to $\mu_2$ most of the time.

It turns out that one needs a generalization of the above approach to
cover some interesting cases. To illustrate this generalization,
consider once again the network of Figure~\ref{figex1}, but assume
that $\lambda> \mu_1 = \mu_2 = \mu$. The vector flow diagram of the
corresponding fluid process is shown in the middle part of Figure \ref
{figex2}. The diagram shows that the fluid process converges to some
point in the set indicated by the two thicker lines: $ \{\bar{h}\}
\times[0, \bar{h}] \cup[0, \bar{h}] \times\{\bar{h}\}$, depending on
the initial condition. While it is true that the rate of the departure
fluid is close to $\mu$ for any point close to that set, it is no
longer the case that the time to reach that limiting set is bounded by
a linear function of the initial distance to the set. For instance, if
the initial state of the fluid process is $(\bar{h}, \bar{h} +
\epsilon
)$ for some arbitrary $\epsilon> 0$, the process takes at least $\bar
{h} / \mu$ to reach the limiting set. To handle this situation, one
considers the set shown in the right-hand part of Figure~\ref{figex2}.
That set has the following two key properties: (1) the departure flow
rate is almost $\mu$ close to that set, and (2) the time to reach the
set is bounded by a linear function of the initial distance to it, as
can be see from the diagram. Thus, as in the previous example, one can
show that the stochastic network has a departure rate close to $\mu$
most of the time.

The main technical contributions of the paper are as follows:
\begin{itemize}
\item A technique for scaling time, space and threshold for finding a
fluid limit for a stochastic network with threshold based ingress
discarding such as in our example;

\item Proof of a fluid limit for stochastic networks with thinned
processes such as $\Lambda(t)$ in Figure~\ref{figex1};

\item Proof of approximation of the rates of the stochastic network by
the rates of the limiting fluid process under the two key properties
indicated in our examples.
\end{itemize}

In the next subsection, we outline the key steps of our analysis. In
Section~\ref{spriorwork} we relate our work to other prior work,
and in Section~\ref{secbigexample} we review an example stochastic
network with ingress discarding. Section~\ref{secpreliminaries}
establishes the notation and initial model description, while
Section~\ref{sfluidanalysis} proves the main results of the article.
In Section~\ref{secswitchexample} we study the fluid model of a
network resembling a $2 \times 2$ network switch and show that the fluid model
has the necessary properties to employ the main results of the article.
Note that Musacchio~\cite{MusacchioPhD} shows that a more general
network with ingress discarding has a fluid model with the necessary
properties. Section~\ref{sconclusion} concludes the paper.\looseness=-1

%As in the works~\cite{Dai95,Bramson,Mandelbaum} and many others, our
%analysis considers a \emph{fluid model} of the stochastic network.
%Roughly, the fluid model of the stochastic network is a system for
%which the discrete dynamics are approximated by continuous dynamics.
%In our context, the fluid model is a system of differential inclusions
%that are dependent on $\bar{h}$.

%Previous work by Dai~\cite{Dai95} and others demonstrate that showing
%stability of the fluid model is sufficient to prove that the
%stochastic network is stable (positive Harris recurrent). Our work is
%similar in that we reduce the analysis of a stochastic network to
%analysis of a fluid model, but what we show is different. We
%demonstrate that if the fluid model is both ``stable'' and achieves
%desired flow rates, then there exists a large enough threshold scale
%factor $n$ such that the long-term average flow rates of the
%stochastic network are arbitrarily close to those of the fluid model.

%s1.1 #&#
\subsection{Proof outline} \label{sproofoutline}
Our goal is to show that the long-term average flow rates of the
stochastic system can be made arbitrarily close to a vector of desired
rates $R$ if the discarding thresholds are made large enough. Moreover,
we want to show that certain properties of the system's fluid model
suffice to reach this conclusion. In this subsection we outline the
arguments detailed in the rest of the paper.

The queuing network we consider has ingress discarding thresholds of
$nh$ in each queue, where $h>0$, and $n>0$ is a threshold scale factor
that is increased to make the thresholds larger. The network is
described by a Markov process $X^n=\{X^n(t), t\geq0\}$ taking values
in the state space $\mathsf{X}$. The superscript emphasizes the
dependence on $n$. The state of the Markov process includes the queue
lengths, remaining service times at each queue and the remaining time
until the next exogenous arrival of each flow $f\in{1,\ldots,F}$. We will
argue that $X^n$ satisfies the strong Markov property.

%To consider systems with ever increasing discarding thresholds, we
%considers systems with thresholds of $nh$ where $n$ is a scaling
%factor.
%The notation $X^{x,n}(\cdot) $ emphasizes the dependence of the
%process on initial state $x$ and parameter $n$.
%T
%wo other processes are used to describe the evolution of the queueing
%network: $\Lambda^{x,n}: [0,\infty) \rightarrow\R^F_+$ that counts
%the number of customers of each flow admitted by the control scheme
%since time $t=0$, and $T^{x,n}: [0,\infty) \rightarrow\R^K_+$ which
%counts the cumulative time (since $t=0$) each of the $K$ queues in the
%network have received service.
%A \emph{fluid model} is defined precisely (later) by a set of
%differential inclusions that are consonant with an intuitive fluid
%analog of the queueing dynamics. A solution to the fluid model is
%called a \emph{fluid model trajectory}. A fluid model trajectory
%consists of three functions, the first of which $\bar{X}: [0,\infty)
%Markov process $X$. ($\mathsf{\bar{X}}$ is similar to state space $
%values in $\R_+$ rather than $\mathbb{Z}_+$.) The other two functions
%that constitute the fluid model trajectory are $\bar{T}: [0,\infty)
%which are analogous to the processes $T^{x,n}$ and $\Lambda^{x,n}$
%respectively.
%For our purposes a ``stable'' fluid model means that all fluid model
%trajectories are drawn to and absorbed by a compact set $\bar{h}\eq
%example of the introduction. The notation $\bar{h}\eq$ emphasizes that
%the set depends on the fluid model discarding threshold, and scales
%multiplicatively with it.

As we discussed in the previous section, we construct fluid limits of
the system by scaling time, space and threshold scale factor in
particular ways that we describe below. These fluid limits converge (in
a sense also described below) to trajectories of a fluid model. \label
{likeoriginal}The fluid model, like the original system, also has
ingress discarding thresholds. However, these thresholds need not equal
$h$, since one of the fluid limits we need to consider can scale space
and threshold at different rates. Therefore when referring to the
system's fluid model, we need to specify~$\bar{h}$, the discarding
thresholds of each queue of the fluid model. (The\vspace*{1pt} queues
of the fluid model have a common threshold $\bar{h}$, just as the
queues of the original system have a common threshold $nh$.) The fluid
model has a state space $\bar{\mathsf{X}}$ similar to that of the
original system, but the queue lengths take values in $\mathbb{R}^+$
rather than $\mathbb{Z}^+$.
%Our goal is to show that if the fluid model satisfies a neis true,
In what follows we adopt the notation that if $\mathcal{S} \subset
\bar{\mathsf{X}}$, then the set $a\mathcal{S}$ ($a \in\mathbb{R}_+$)
denotes a ``scaled'' set such that $ \bar{x} \in a\mathcal{S}$ iff $
\bar{x}/a \in\mathcal{S}$. Also let $\llVert \bar{x} \rrVert_{
\mathcal{S}} =
\inf_{e \in\mathcal{S} } \|\bar{x}-e\|$ denote the distance between
$\bar
{x}$ and the set $\mathcal{S}$.

Our goal is to show that if there exists a closed, bounded set \label
{firsteq} $\eq\subset\bar{\mathsf{X}}$ and $t_0 \in\mathbb{R}^+$
such that conditions (C1) and (C2) below hold, then there exists a large
enough $n$ such that the stochastic network achieves long-term average
departure rates arbitrarily close to $R$.
\label{conditions}Conditions (C1) and (C2) are as follows:\looseness=-1
\begin{longlist}[(C2)]
\item[(C1)] All trajectories of the fluid model with ingress
discarding thresholds $\bar{h}$ and initial condition $\bar
{X}(0)=\bar
{x}$ are absorbed by a set $\bar{h}\eq$ in a time not more than $t_0
\llVert \bar{x} \rrVert_{ \bar{h} \eq} $;
\item[(C2)] If $\bar{h}>0$, the instantaneous departure rates of
the fluid model while its state is in the set $\bar{h} \eq$ are equal
to the vector of desired rates $R$.
\end{longlist}
%
%}
%Note that (C1) requires that $\bar{h} \eq$ be
\label{absorbing}Note that (C1) requires that $\bar{h}\eq$ be an
absorbing set of the fluid model with thresholds $\bar{h}$. For
example, one can show that a minimal absorbing set of the fluid model
in many cases would be, roughly, the set of states such that at least
one of each flow's set of ``bottleneck'' queues is at it's discarding
threshold, and servers with a utilization below 1 have empty queues.
(By ``bottleneck queue,'' we mean a queue whose service constrains a
flow's rate in the fluid model.) However, such a construction might not
be sufficient to satisfy (C1), particularly when flows do not have unique
``bottlenecks.''
Recall that in the \hyperref[sintro]{Introduction} we studied an example with two
serially-connected queues with the same service rate.
This is an example in which $\bar{h} \eq$ needs to be made larger than
the minimal absorbing set in order to satisfy (C1).
To see this note that even though the two line segments in the middle
panel of Figure~\ref{figex2} constitute an absorbing set for the fluid
model, if we defined $\eq$ so that $\bar{h}\eq$ is equal
to these two line segments (by making $\eq= \{ 1 \} \times[0,1] \cup
[0,1] \times\{ 1\} $), condition (C1) would not be met.
By defining $\eq$ in such a way as to make $\bar{h} \eq$ have the shape
indicated by the shaded area of the right panel of Figure \ref
{figex2}, the time it takes trajectories of the fluid model to reach
$\bar{h} \eq$ can be upper bounded by an amount proportional to the
distance of the starting point of the trajectory from $\bar{h} \eq$,
thus satisfying (C1).

%Our goal is to show that if the fluid model of the system meets
%certain conditions, then there exists a large enough $n$ such that the
%stochastic network achieves long-term average departure rates
%arbitrarily close to a vector of desired rates $R$. Let $\eq\subset
%mathsf{X}$ be an arbitrary closed, compact set

%%We now state these necessary conditions. Like the original network,
%the fluid model has ingress discarding thresholds. Consider the fluid
%model with unit thresholds and let $\eq\subset mathsf{X}$ be a
%closed, compact set. % that contains the invariant set of the fluid
%model.
%%To state these conditions, we adopt the following notation. Let $\eq$
%be a set in the state space of the fluid model that contains the
%invariant set of the fuid
%The notation $\bar{h}\eq$, $\bar{h} \geq0$ denotes a scaled set such
%that $ \bar{x} \in\bar{h}\eq$ iff $ \bar{x}/\bar{h} \in\eq$. Also
%let $\n{x}{\bar{h} \eq} = \inf_{e \in\bar{h}\eq} ||x-e||$ denote the
%distance between a state $x$ and the set $\bar{h}\eq$. The conditions
%are:
% \begin{itemize}
% \item[\textbf{(C1)}] All trajectories of the fluid model with ingress
%discarding thresholds $\bar{h}$ and initial condition $\bar{X}(0)=
%$t_0 \n{\bar{x}}{\bar{h} \eq} $ where $t_0$ is proportionality
%constant;
% \item[\textbf{(C2)}] The instantaneous departure rates of the fluid
%model while in states in the set $\bar{h} \eq$ %(which depend on the
%first derivative of the $\bar{T}(\cdot)$ process)
% are equal to a vector of desired rates $R \in\R_+^F$, provided $
% \end{itemize}
%Note that (C1) requires that $\bar{h}\eq$ contain the invariant set

The proof depends on two main steps: %In both steps we consider the
%network's Markov process description $X(\cdot)$ scaled in time and
%space by the size of the threshold scale factor $n$. This scaled
%process $Y(\cdot)$ satisfies $Y(t) = X(nt)/n$ for each $t\in[0,
%
\begin{longlist}[(ii)]
\item[(i)] The expected flow rates associated with the process
$X^n(\cdot)$, over a finite time interval of length $nt_0$, and for
initial conditions near a set $ n h \eq$, can be made to be arbitrarily
close to $R$ with a sufficiently large threshold scaling factor $n$.

\item[(ii)] The excursions of the process $X^n(\cdot)$ away from $n h
\eq$ become relatively shorter with larger threshold scaling factor
$n$. More precisely, the first hitting time that occurs $n t_0$ after
having started in a neighborhood of the set $n h \eq$, can be made to
be arbitrarily close to $n t_0$.
\end{longlist}

In both steps we make use of the fact that a fluid limit of the process
$X^n(\cdot)$ converges to a trajectory of the fluid model, but the
different objectives of the two steps require us to use different fluid
limit scalings. In the first step we consider a sequence of (initial
condition, scale factor) pairs $\{(x_j,n_j)\}$. To emphasize the
dependence on initial condition and threshold scale factor we write $
X^{\mathbf{x}_j}(\cdot)$, where the superscript ${{\mathbf {x}}_{j}}
\triangleq(x_j,n_j)$. We require that the sequence has the properties
that $x_j / n_j$ is no more than a distance $\zeta<1$ away from the set
$h\eq$, and $n_j \rightarrow\infty$. Otherwise, the sequence is
arbitrary. We call such a sequence a \textit{near fluid limit}
sequence. (Equivalently, the near fluid limit condition has $\llVert
x_j \rrVert_{ n_j h \eq} < n_j \zeta$ and $n_j \rightarrow\infty$. In
general it is often more intuitive to consider the distance of $X/n$
from the set $h \eq$ than to consider the distance of $X$ from $n h
\eq$, so we will use whichever construction is more convenient or
intuitive for the context.) We demonstrate that the sequence of scaled
processes $\{\frac{1}{n_j} X^{\mathbf{x}_j}(n_j \cdot)\}$ converges
along\vspace*{-1pt} a subsequence, uniformly over compact time
intervals, to a fluid model trajectory $\bar{X}(\cdot)$. The result
largely follows from the fact that the process describing the
cumulative time each server in the network is busy is Lipschitz
continuous, and a sequence of Lipschitz continuous functions on a
compact set converges along a subsequence. Consequently, the
convergence to a fluid trajectory only\vspace*{1pt} holds on a finite time interval.
The thresholds of the fluid model that $\bar{X}(\cdot)$ satisfies are
of size $\bar{h}=h$. This is because we scale both space and threshold
by the same amount in this fluid limit, so the two scalings cancel out.
Moreover, the restrictions we put on the near fluid limit sequence
ensure that the initial condition of the fluid model trajectory $\bar
{X}(\cdot)$ is within\vspace*{1pt} a distance of $\zeta$ of $\bar{h}\eq$. Thus, the
fluid model trajectory $\bar{X}(t)$ hits $\bar{h}\eq$ quickly [in not
more than time $\zeta t_0$ by (C1)] and then achieves flow rates of $R$
[by condition (C2)].

At this point, we have only shown convergence along a subsequence to a
fluid trajectory with some desired properties. We need to show
convergence along the original near fluid limit sequence in order to
eventually make conclusions about the stochastic network. To that end,
consider a functional $\frak{F}$ that extracts the difference between
the actual flow throughput and the desired flow throughput over a
compact time interval $[\zeta t_0, t_0]$ (in time scaled by $n$).
Since $\bar{X}(t)$ hits $\bar{h} \eq$ by time $\zeta t_0$, the flow
rates are equal to the desired rates over $[ \zeta t_0, t_0 ]$.
Consequently, $\frak{F} \circ\bar{X}=0$. This in turn allows us to
argue that $\{ \frak{F} \circ\frac{1}{n} X^{\mathbf{x}_j}( n
\cdot) \}$
converges to 0 along a subsequence. Since every near fluid limit
sequence of processes (with the functional applied to them) converges
along a subsequence to 0 in this way, it must be that every near fluid
limit sequence also converges to 0 in this way. This fact allows us to
show that the flow rates of the process $\frac{1}{n}X(n \cdot)$ can be
made arbitrarily close the desired rates, for a finite time period,
from any scaled initial condition $x/n$ near~$h\eq$, provided that $n$
is sufficiently large. In the detailed proof the functionals we
consider act on the Markov state trajectory combined with the
trajectories of some other associated processes such as the cumulative
service time process. \label{factthat}The fact that $\zeta<1$ was
chosen otherwise arbitrarily is important because it allows us to later
make $\zeta$ small so that the desired rates are achieved over most of
the interval $[0, t_0]$ (in scaled time).

In the second step, we again consider a sequence of (initial condition,
scale factor) pairs $\{(x_j,n_j)\}$.
This sequence must satisfy the properties that the distance between
$x_j/n$ and $h\eq$ is more than a constant $\zeta$ for each $j$, and
that $ \llVert x_j \rrVert_{ n_j h \eq} =n_j \llVert x_j/n_j
\rrVert_{ h\eq}\rightarrow\infty$.
%%(The later condition is equivalent to requiring that $\n{x_j}{n_j h
%
Otherwise, the sequence is arbitrary.
We call such a sequence a \textit{far fluid limit} sequence.
We show that the sequence of scaled processes $\{X^{\mathbf{x}_j} (
\llVert x_j \rrVert_{ n_j h\eq} \cdot)/ \llVert x_j \rrVert_{
n_j h \eq} \}$ converges along a
subsequence of any far fluid limit sequence, uniformly over compact
time intervals, to a fluid model trajectory $\bar{X}(\cdot)$ satisfying
a fluid model with discarding thresholds $\bar{h}$. The scaled
threshold sequence of the fluid limit is $\{ n_j h / (n_j \llVert
x_j/n_j \rrVert_{ h \eq} ) \}$, so the choice of sequence and
convergent subsequence
determines a value for $\bar{h}$ that satisfies $\bar{h}\in[0,h\zeta
^{-1}]$. Also the scaling of the far fluid limit sequence ensures that
the initial condition of the fluid trajectory have an initial condition
that is unit distance from $\bar{h}\eq$.
This fact along with our starting assumption (C1), ensure that $\llVert
\bar{X}(t_0) \rrVert_{ \bar{h} \eq} = 0$.
The preceding two facts allow us to argue that the sequence $\{
X^{\mathbf{x}_j} ( \llVert x_j \rrVert_{ n_j h \eq} t_0)/ \llVert
x_j \rrVert_{ n_j h \eq} \}$ has a distance
from $n_j h \eq$ that converges to 0 along a subsequence.
Moreover since any far fluid limit sequence has a subsequence that
converges to 0 in this sense, it must be that this convergence property
holds for any far fluid limit sequence.

This fact is the basis for constructing an argument that
\[
\E\bigl\llVert X^{{\mathbf{x}}_{}} \bigl( t_0 \llVert x
\rrVert_{ n h
\eq} \bigr) \bigr\rrVert_{ n h \eq} \leq\delta\llVert x
\rrVert_{ n h \eq}
\]
for any $\delta>0$ provided that threshold scale factor $n$ is
sufficiently large and $\llVert x/n \rrVert_{ h\eq}>\zeta$
(equivalently $\llVert x \rrVert_{ nh\eq} > n \zeta$). This
relation serves as a Lyapunov function which allows
the construction of an argument about the recurrence time of the scaled
process $X/n$ to a neighborhood with distance $\zeta$ of $h\eq$, and
this in turn allows us to conclude (ii) above.

This recurrence time argument is adapted from~\cite{MeynTweedie} while
the overall argument we make with the far fluid limit sequence
parallels~\cite{Dai95}. The main difference between our far fluid limit
argument and that of~\cite{Dai95} is that in~\cite{Dai95} the fluid
model and stochastic network are drawn to the origin and neighborhood
of the origin, respectively, whereas in our model the system is
attracted to a set of states.

% In this overview of the proof technique, we of course omitted certain
%details in order to provide a concise summary. For instance the
%convergence to a fluid trajectory, and other outline above, make use
%of cumulative service time process, in addition to the Markov state
%process $X$.

%s1.2 #&#
\subsection{Relation to prior work} \label{spriorwork}

Our fluid limit proof techniques borrow heavily from work by Dai \cite
{Dai95}. Dai shows that for networks without discarding, stability of a
corresponding fluid model implies positive Harris recurrence of the
stochastic network. In our work we use the fluid model not only to show
positive Harris recurrence of the stochastic network, but also to find
its long term average flow rates. Specifically, we use two fluid
limits: the far fluid limit and the near fluid limit that correspond to
different sequences of initial conditions and threshold pairs.

Dai's proof considers a sequence of initial states $\{x\}$ of the
Markov process describing the network, with $|x| \rightarrow\infty$,
and then obtains a fluid limit by scaling time and space by $|x|$. Dai
uses this result to construct a Lyapunov function to show that the
expected state of the system contracts, for initial states far enough
from the origin. Our far fluid limit analysis parallels this, but with
the difference that our analysis focuses on the distance of the state
from a set of states $h \eq$ rather than the distance from the origin.
Also, because we are interested in showing the existence of a
sufficiently large threshold scaling factor $n$, for both the near and
far fluid limits, we consider a sequence of initial condition threshold
pairs $\{x,n\}$ to obtain our results rather than just a sequence of
initial conditions as in~\cite{Dai95}.

Our fluid limit technique is also very similar to that found in work by
Bramson~\cite{Bramson}. In much the way we do, Bramson takes the fluid
limit using a sequence of pairs, one being the initial condition and
the other being a time scaling factor of both space and time. However,
our results do not follow immediately from the results of Bramson
because we require that the fluid model be drawn toward a set $\bar
{h}\eq$ rather than just to the origin.

Another body of work uses fluid limits to show rate stability rather
than showing that the system state converges to an invariant
distribution, or more precisely that the system is positive Harris
recurrent. Rate stability means that the long-term average departures
match the long-term average arrivals. It is a weaker concept than
positive Harris recurrence because a system can be rate stable while
internally the average queue lengths grow unbounded or at least fail to
converge to an invariant distribution. For a treatment see \cite
{ElTaha99}, and examples of its application include~\cite{Chen95} and
\cite{Dai00}. The rate stability framework is not sufficient for our
objectives because in order to show that our control policy achieves
flow rates close to those predicted by a fluid model, we need to show
that the vector of queue lengths settles to an invariant distribution
concentrated near a particular set of lengths, as illustrated in the
example of the \hyperref[sintro]{Introduction}.

Another closely related work to ours is by Mandelbaum, Massey and
Reiman~\cite{Mandelbaum}. In~\cite{Mandelbaum}, the authors study the
fluid limit of a queueing network with state dependent routing, where
the function describing the arrivals to each queue can scale with $n$
and or $\sqrt{n}$, in a manner similar to the scaling of our
thresholds. The authors prove a functional strong law of large numbers
and a functional central limit theorem in the context of their model.
However, the authors assume that the network is driven by Poisson
processes, rather than just the renewal assumption that we make. An
earlier work by Konstantopoulos, Papadakis and Walrand derives a
functional strong law of large numbers and a functional central limit
theorem for networks with state dependent service rates~\cite{Takis}.

There are also several other works that use reflected Brownian motion
models to study queueing networks with blocking \cite
{DaiHarrison,Harrisonbook,harrisonrbm}. Typically the objective of
most such investigations is to approximate the distribution of the
queue occupancy with a diffusion approximation. In contrast with those
works, our objective is to show almost sure convergence using a strong
law of large numbers scaling.

%s1.3 #&#
\subsection{Example network}
\label{secbigexample}

In this subsection, we introduce an example that motivates the theory
developed in this paper. The example will illustrate two important
phenomena---that the long-term rates of the stochastic system get
closer to those of a corresponding fluid model when discarding
thresholds are raised, and that when there are not unique bottlenecks,
the vector of queue depths is not attracted to a unique
equilibrium point.

%f3 #&#
%
\begin{figure}

\includegraphics{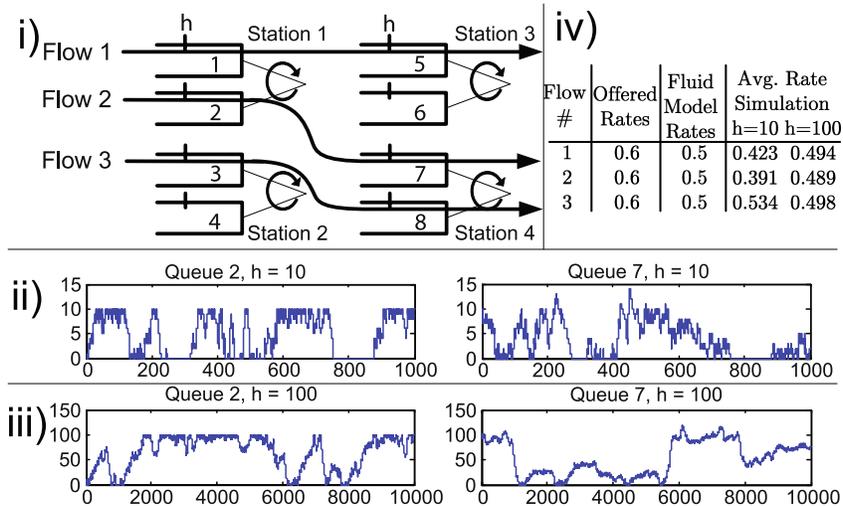}

\caption{\textup{(i)} An illustration of the example network. \textup{(ii)} The queue
lengths of queue 2 and 7 when the discarding threshold parameter $h$
set to $10$. \textup{(iii)} The lengths of queue 2 and 7 when the discarding
threshold parameter $h$ is set to $100$. \textup{(iv)} The average flow rates of
each of the three flows for both the $h=10$ and $h=100$ simulated
sample paths.}
\label{figexample2}
\end{figure}

Our example is illustrated in Figure~\ref{figexample2}. The example is
analogous to a two-input and two-output switch. Two flows enter the
network at station 1, the first input of the switch, and a third flow
enters the network at station 2. We concentrate on flow~2, which shares
stations 1 and 4 with flows 1 and 3, respectively. All stations are
served at rate 1, have round-robin service with equal weighting to all
queues, and have service times that are exponentially distributed. The
arrival rate of each flow is 0.6, with Pareto inter-arrival
distributions given by
\[
\Pb\bigl(\xi_f(j)> s\bigr) = \frac{1}{(0.6 s+1)^{2} },\qquad f\in\{1,2,3\},
s \geq0, %\label{eqpareto}
\]
where $\xi_f(j)$ is the inter-arrival time preceding the $j$th arrival.
We choose the Pareto distribution for this example to emphasize that we
are interested in networks whose inter-arrival and service times are
not necessarily memoryless.

We consider the behavior of the network's fluid model. Since stations 2
and 3 have a capacity of 1 and each carry one flow with an offered rate
of 0.6, the queues of these stations should never fill. Stations 1 and
4 each carry 2 flows that offer a load of 0.6 (before considering
discarding). The fluid model of the station's round robin service is
that each station serves both of its queues at rate 0.5 as long as both
flows are offering enough customers to be served at this rate.
Consequently, when flow 1's queue at station 1 is filled below
threshold, this queue grows at a rate of 0.1. However, if flow 1's
queue at station 1 ever went above its threshold, ingress discarding
would commence, and the queue would immediately decrease. Therefore, it
must be that this queue grows to its threshold, stays at this level and
then flow 1's ``thinned'' or post-discarding arrival process is of rate 0.5.

Similar reasoning shows that flow 3's queue at station 4 behaves in
this way, and also that one of flow 2's queues must also reach the
threshold and ``stick'' there. These steps allow us to conclude that
after some time, all three flows should have rates of 0.5 in the
system's fluid model. (We will verify this carefully in Section~\ref
{secswitchexample}.)

Figure~\ref{figexample2} shows the simulated trajectories of flow 2's
queues at both bottleneck stations in the stochastic network. In the
$h=10$ case, the simulation shows that queues 2 and queue 7, which both
serve flow 2, are empty for over 100 time units around time 800. This
empty period is significant because when flow 2's queues are empty,
flow 2 misses opportunities to have its customers served by the
bottleneck stations. Indeed the table included in Figure~\ref
{figexample2} shows the average rate, averaged over the last $80\%$ of
the simulation time to reduce some of the initial transient effect, is
$0.391$. This is substantially below the rate of $0.5$ predicted by the
fluid model. Most likely, a string of long interarrival times of flow
$2$, caused the queues at the bottleneck stations to starve.

Raising the thresholds should reduce starvation, because larger
thresholds would provide the bottleneck queues a larger backlog to
smooth over fluctuations in the arrival and service processes. To test
that intuition, we simulate the network with discarding thresholds of
$h=100$. Figure~\ref{figexample2} shows the trajectories of flow 2's
queues for the increased threshold. We note that neither queue spends
all of the time filled to its threshold, but instead at most times at
least one of the queues is near its threshold. For instance, at the
beginning of the simulation, queue 7 (the second bottleneck) is
chattering near the threshold while queue 2 (the first bottleneck) is
below threshold. At some time before the 2000 second mark, the two
queues switch these roles, and around the 6000 second mark the queues
switch these roles again. We also note that flow 2 achieves an average
rate of $0.489$, which is much closer to the rate of $0.5$ predicted by
the fluid model.

%s2 #&#
\section{Preliminaries}
\label{secpreliminaries}
Customers of a given flow $f \in\{1,\ldots,F\}$ follow the same fixed
sequence of distinct stations. The service times are independent. Each
flow $f$ has a weight $w_f$ and each station $i\in\{1,\ldots,d\}$ is
equipped with per-flow queues and serves each flow in proportion to its
weight using a weighted round robin or a similar queueing discipline.
In addition to the notion of flow, each customer also has a class $k\in
\{1,\ldots,K\}$ that is indicative of both the customer's flow and the
station $s(k)$ it is located. Thus the class of a flow $f$ customer
changes as the customer progresses from station to station, but with
the restriction that a flow $f$ customer must always have a class in
the set $K(f)$. Conversely, each class $k$ is associated with one
and only one flow $\ff{k}$. We also adopt the numbering convention that
flow $f$ customers enter the network as class $k=f$, and thus $f\in
K(f)$. The constituency matrix $C \in\{0,1\}^{d \times K}$ records
which classes are served in each station: $C_{ik}=1$ if class $k$ is
served at station $i$, otherwise $C_{ik}=0$. A customer of class $k$
who completes service becomes a customer of class $l$ if $P_{kl}=1$.
Thus $P \in\{0,1\}^{K\times K}$ is a binary incidence matrix with each
row containing at most one~1. Because flows follow loop-free paths, $P$
is nilpotent.

The exogenous arrivals to the network for flow $f$ are described by a
renewal process $E_f(\cdot)$ for which the interarrival times $\{\xi
_f(j),j \geq1 \}$ are i.i.d. and $\alpha_f$ is the mean arrival rate. Thus,
\[
E_f(t) = \max\bigl\{r\dvtx  U_f(0)+\xi_f(1) +
\cdots+ \xi_f(r-1) \leq t \bigr\},\qquad t\geq0,
\]
where $U_f(t) \in\R_+$ is the time after time $t$ until the next flow
$f$ customer arrives at the network ingress. We also need to assume
that inter-arrival times are unbounded and spread-out. More precisely,
we assume that for each $k\in\{0,\ldots,F\}$, there exists an
integer $j_k$ and some function $p_k(x)\geq0$ on $\R^+$ with $\int
_0^\infty p_k(x) \,dx >0$, such that
%
%e1 #&#
%
\begin{equation}
\label{equnbounded} \Pb\bigl[ \xi_k(1) \geq x \bigr] > 0
\qquad\mbox{for any } x>0
\end{equation}
and
%
%e2 #&#
%
\begin{equation}
\label{eqspreadout} \Pb\Biggl[a \leq\sum
_{i=1}^{j_k} \xi_k(i) \leq b \Biggr] \geq
\int_a^b p_k(x) \,dx \qquad\mbox{for any }
0\leq a \leq b.
\end{equation}

The service times $\{\eta_k(j),j \geq1 \}$ of each class $k$ are also
i.i.d. and have mean $m_k=\mu_k^{-1}$, where $\mu_k$ is the mean
service rate. We also define the $K \times K$ diagonal matrix $M$ whose
$k$th diagonal entry is $m_k$. The quantity $V_k(t) \in\R_+$ denotes
the remaining service time of the class $k$ customer in service, if
there is one at time $t$, otherwise $V_k(t)=0$. We define a service
process $S^{\mathbf{x}}_k(\cdot)$ as
\[
S_k(t) \triangleq\max\bigl\{j\dvtx  \tilde{V}_k(0) +
\eta_k(1) +\cdots+\eta_k(j-1) \leq t\bigr\},\qquad t\geq0,
\]
where $\tilde{V}_k(0) = V_k(0)$ if $V_k(0)>0$; otherwise $\tilde
{V}_k(0)= \eta_k(0)$ is a fresh service time with the same distribution
as $\eta_k(1)$ and independent of all other service times.

In principle, our assumption that the service times are independent
does not allow for service times that depend on a packet's size (taking
``packets'' to be ``customers''). Dependence on packet size would make
the service times of stations dependent on each other. To model this
explicitly would require a much more complicated model. However we
believe that our results in this work would still hold if this
assumption were relaxed.

We define the following right-continuous processes: $A\dvtx[0,\infty)
\rightarrow\mathbb{Z}_+^K$ counts the arrivals to each class $k$ since
time $t=0$; $D\dvtx[0,\infty) \rightarrow\mathbb{Z}_+^K$ counts the
departures of each class; $\Lambda\dvtx[0,\infty) \rightarrow
\mathbb
{Z}_+^F$, counts the exogenous arrivals of each flow that make it past
the discarding point (``thinned'' exogenous arrivals); $Q\dvtx
[0,\infty
) \rightarrow\mathbb{Z}_+^K$ is the vector process of queue depths;
$T\dvtx[0,\infty) \rightarrow\mathbb{R}_+^K$ counts the total time each
class $k$ has been served since $t=0$; and $I\dvtx[0,\infty)
\rightarrow
\mathbb{R}_+^d$ counts the total time each server has been idle since
$t=0$. For each $t\geq0$, these processes satisfy the following relations:
%
%e3 #&#
%e4 #&#
%e5 #&#
%e6 #&#
%e7 #&#
%e8 #&#
%e9 #&#
%
\begin{eqnarray}
\label{eqarrivaldef}
&\displaystyle A(t) = P^T D(t) + \Lambda(t);&
\\
\label{eqqevolution}
&\displaystyle  Q(t) = Q(0) + A(t) - D(t);&
\\
\label{eqQbig0}
&\displaystyle Q(t) \geq0;&
\\
\label{eqTnd}
&\displaystyle T_k(t) \mbox{ is nondecreasing and }
T_k(0)=0\qquad \mbox{for $k=1\lldots K$};&
\\
\label{eqInd}
&\displaystyle I_i(t) = t - C_iT(t) \mbox{ is
nondecreasing and } I_i(0) =0\qquad \mbox{for $i=1\lldots d$};&
\\
\label{eqIdleLimit}
&\displaystyle \int_0^{\infty} \bigl( C Q(t)
\bigr) \,dI(t) = 0;&
\\
\label{eqdepartures}
&\displaystyle  D_k(t) = S_k\bigl(T_k(t)
\bigr) \qquad\mbox{for $k=1\lldots K$}.&
\end{eqnarray}
%
%The expression $C_i$ in (\ref{eqInd}) is the $i$th row of the
%constituency matrix $C$.
Relations (\ref{eqarrivaldef})--(\ref{eqQbig0}) describe the
relations between the arrival, departure and queue length processes.
Statements (\ref{eqTnd})--(\ref{eqIdleLimit}) describe basic
restrictions on the cumulative service time and idle time processes,
with relation (\ref{eqIdleLimit}) reflecting an assumption that each
station is work conserving. Equation (\ref{eqdepartures}) reflects
that departures of class $k$ are determined by the composition of the
service time counting process $S_k(\cdot)$ and the process $T(\cdot)$.

The ingress discarding scheme drops arriving customers of flow $f$ as
they arrive whenever any queue in the set $K(f)$ exceeds a high
threshold $nh$. Recall that $n$ is the threshold scaling factor which
we will adjust in our analysis. \label{Conversely}Conversely, when all
of the queues in $K(f)$ are below a lower threshold $nh-o(n)$, flow $f$
customers are permitted to enter the network. Note the lower threshold
could be set to be the same as the upper threshold, but in some
practical applications it might be beneficial to have different
thresholds so that the switching between admitting and discarding is
less frequent. Thus we permit this difference between upper and lower
thresholds to be any function $o(n)$ that satisfies $o(n)/n \rightarrow
0$ and $o(n)\geq0$. For instance any nonnegative constant may be used.
Between these thresholds, the system has hysteresis behavior, and we
define this behavior as follows. A process $H_k\dvtx[0,\infty)
\rightarrow\{0,1\}$ keeps track of whether discarding has been
``turned-on'' by each class $k$ queue. If $Q_k(t) \geq nh$, then
$H_k(t) =1$, and if $Q_k(t)\leq nh-o(h)$, then $H_k(t)=0$. For all $t$
such that $Q_k(t)\in(nh-o(h),nh)$, the evolution of $H_k$ is
determined by the following rules:
\begin{itemize}
\item If $H_k(t)=0 $, then let $t_s = \min\{\tau\geq0\dvtx  Q_k(t+\tau
)\geq nh\}$ [note that $t_s$ is well defined because $Q_k(\cdot)$ is right
continuous]. $H_k(t+ \tau)=0$ for $\tau\in[0,t_s)$ and $H_k(t_s) = 1$.
\item If $H_k(t)=1 $, then let $t_s = \min\{\tau\geq0\dvtx  Q_k(t+\tau) <
nh-o(n) \}$. $H_k(t+\tau)=1$ for $\tau\in[0,t_s)$ and $H_k(t_s) = 0$.
\end{itemize}

The flow $f$ customers that are allowed into the network beyond the
discarding point depends on all the processes $H_k(\cdot)$ as
%
%e10 #&#
%
\begin{equation}
\label{eqLambda} \Lambda_f(t) = \sum_{j=1}^{E_f(t)}
\prod_{K(f)} \bigl(1- H_k(
\tau_j-)\bigr),
\end{equation}
where $\tau_j = U_f(0) + \sum_{m=1}^{j-1} \xi_f(m) $ is the time of
the $j$th
arrival to the discarding point. Here the dependence on $H_k(\tau_j-)
\triangleq\lim_{t \uparrow\tau_j} H_k(t)$ rather than $H_k(\tau_j)$
is to avoid problems with causality. For instance a customer arrival
that triggers discarding should not be discarded; otherwise the
customer will never arrive to the system, and paradoxically the
discarding will never turn on. Our modeling choice allows such a
customer to enter, thus triggering discarding, which will discard
future customers.

The queueing discipline of a station $i$ serves each flow in proportion
to the flow weights over long time intervals. More precisely, for some
constant $c>0$ and all $\tau>0$,
%
%e11 #&#
%
\begin{equation}
\label{eqRRdrainf} \frac{D_k(t,t+\tau)}{w_{\ff{k}}} \geq\frac
{D_l(t,t+\tau)}{w_{\ff{l}}} - c \qquad\mbox{whenever
$Q_k(s)>0$ $\forall s\in[t,t+\tau]$}\hspace*{-20pt}
\end{equation}
for all $k,l \in C(i) \triangleq\{k'\dvtx C_{ik'}=1\}$, where $D_k(t,t+\tau
) \triangleq D_k(t+\tau) - D_k(t)$.

We furthermore assume that only the customer at the head of line of
each queue may be served, and that the instantaneous service rate of
any queue is a function of the current state. That is,
%$
$\dot{T}_k(t) = f(X(t))$ for some function $f(\cdot)$
%$
where $X(t)=[Q(t);U(t);V(t);H(t)]$.

The evolution of the queuing system depends on the particular queuing
discipline. Moreover, some queueing disciplines require additional
state variables. For instance, a weighted round robin scheduler visits
the queues in a cyclic order, serving any customers at the head of the
line. The order should be chosen so that in each cycle the number of
visits of each queue is proportional to the flow weights (which is
possible if the weights are rational multiples of each other). Other
queueing disciplines could be considered as well, though these
disciplines may need additional state variables. For instance, deficit
round robin (DRR) requires counters for each class~\cite{DRR}. Also,
DRR ensures that the service times given to each class are proportional
rather than the number of customers served.\vadjust{\goodbreak} Therefore, DRR satisfies a
criterion similar to (\ref{eqRRdrainf}) except that $D(\cdot)$ is
replaced by $T(\cdot)$. However, since the service times are unbounded,
the criterion holds only in the limit $\tau\rightarrow\infty$, almost
surely. Other disciplines require yet more complex state descriptions.
For instance, weighted fair queueing (WFQ) keeps track of each
customer's ``virtual finish time''---the time they would have departed
if the service discipline were weighted processor sharing and no more
customers were to arrive~\cite{Bertsekas}. To keep the presentation
simple, we assume that the additional state variables required by the
queueing discipline are described by a bounded vector in $\mathbb
{Z}_+^d$. We append this to the $H$ portion of the state description.
Treatment of queueing disciplines that require more elaborate state
descriptions requires some modification to the statement and proof of
Theorem~\ref{thsubseqconv}.

\subsection{State description} \label{sstate}

The dynamics of the queueing network are described by the Markov
process $X=\{X(t),t\geq0\}$. The state description contains the
queue lengths $Q(t) \in\R_+^K$ of all the $K$ queues in the network,
as well as the residual arrival and service times $U(t) \in\R_+^F$ and
$V(t) \in R_+^K$, respectively. Recall that $U(t)$ and $V(t)$ are
defined to be right-continuous. Finally the state description includes
the state of the discarding hysteresis and any state variables used by
the queueing discipline as described above. We assume that $H(t) \in\{
0,1\}^K \times\mathbb{Z}_+^d$. Thus the full state description is
\[
X(t) = \bigl[ Q(t);U(t);V(t);H(t) \bigr].
\]
Let $\mathsf{X} \subset\mathbb{Z}^K \times\mathbb{R}_+^{F+K}\times
\{
0,1\}^K \times\mathbb{Z}_+^d$ be the set of all states $X$ can take.
A~fixed threshold scaling factor $n$, an initial condition $x=X(0)\in
\mathsf{X}$ is sufficient to specify the statistics of the future
evolution of the system.

We claim that the process $X$ satisfies the strong Markov property, by
the same argument given by Dai~\cite{Dai95}. In turn, Dai's argument
followed from Kaspi and Mandelbaum~\cite{kaspimandelbaum}. Without
repeating all the details of the argument, the basic idea is that
$X(\cdot)$ is a piecewise deterministic Markov (PDM) process---behaving
deterministically between the generation of ``fresh'' inter-arrival or
service time. Davis shows that a PDM process whose expected number of
jumps on $[0,t]$ is finite for each $t$ is strong Markov~\cite{davis}.
As we assume that the inter-arrival and service times have a positive
and finite mean, the expected number of jumps of $X(\cdot)$ in any
closed time interval is finite. Therefore $X(\cdot)$ has the strong
Markov property.

The fluid\vspace*{1pt} model, whose defining equations will be given in Theorem \ref
{thsubseqconv}, takes values in the state space $\bar{\mathsf{X}}
\subset\mathbb{R}_+^{F+3K+d}$ since integer valued states of the
original system correspond to real valued states of the fluid model.

%s3 #&#
\section{Fluid limit analysis}
\label{sfluidanalysis}

In this and subsequent sections, we use the superscript ${{\mathbf
{x}}_{}}\triangleq(x,n)$ to denote the dependence on initial state $x$ and
threshold scaling factor $n$.\vadjust{\goodbreak} As we discussed earlier, we use two
different fluid limits in our analysis: the near and far fluid limits
that study behavior of the stochastic network for scaled initial
conditions near and far from $h\eq$, respectively. Recall that $\eq
\subset\bar{\mathsf{X}}$ is a closed and bounded set. Also recall
$h\eq= \{x\dvtx  x/h \in\eq\}$. At this point we make no further
assumptions on $ \eq$, but eventually $ \eq$ will have to be chosen so
that $\bar{h} \eq$ is an absorbing set of the fluid model with
thresholds $\bar{h}$ to apply our final results.

% Much of the analyses of the near and far fluid limits are similar.
%For convenience, we define a scaling function $\s{}$ that allows us to
%consider both
%fluid limits at once. The function $\s{}$ is defined by\index{$\langle
% \[ \s{}_\const\triangleq
% \cases{
% n, & \mbox{if $\n{x/n}{h\eq} \leq\const$} \\
% n\n{x/n}{h\eq}, & \mbox{if $ \n{x/n}{h\eq} > \const$}
% }
% \]
%where the $\const$ parameter is an arbitrary positive constant. We
%will omit the $\const$ subscript when either the choice of $\const$
%does not matter, or when its value is clear from the context.
For notational convenience we also define an augmented state vector
process
\[
\frak{X}^{{\mathbf{x}}_{}}{(\cdot)} \triangleq\bigl[X^{\mathbf
{x}}(\cdot
);T^{\mathbf{x}}(\cdot); \Lambda^{\mathbf{x}}(\cdot);nh\bigr],
\]
which contains all the functions that we want to show converging in
both kinds of fluid limit.

%Before proving Theorem 1, we need a few preliminary lemmas.

In this section, we state Theorem~\ref{thsubseqconv} which shows
convergence to a fluid model trajectory along a fluid limit. The
convergence of the trajectory is uniformly on compact sets. More
precisely, we say that $f_j(t) \rightarrow f(t)$ \textit{uniformly on
compact sets} (u.o.c.) if for each $t\geq0$,
\[
\lim_{j\rightarrow\infty} \sup_{0\leq s \leq t} \bigl|f_j(s) - f(s)\bigr| = 0.
\]
We also use the notation $\dot{f}(t) = \frac{d}{dt}f(t)$ where such a
derivative exists. If a function $f(\cdot)$ is differentiable at $t$,
we say that $t$ is a \textit{regular point}.

The proof, along with four lemmas used in the proof, are given in the
\hyperref[app]{Appendix}. One of these lemmas, Lemma~\ref{lemthin}, is a new result
showing that the thinned arrival process converges u.o.c. to the fluid
limit. In Section~\ref{sproofoutline} we previewed the two types of
fluid limits, which we call the ``near'' and ``far'' fluid limits, that
we will use in our analysis. In both types of fluid limits, time and
space is scaled by a factor that increases. In the development that
follows, that scale factor for time and space is represented by the
notation $a_j$. Later on, we will make specific assumptions about $a_j$
that correspond to either the near or far fluid limit. Bramson \cite
{Bramson} takes a similar approach to defining the fluid limit.
Both types of fluid limit scale the threshold no faster than time and
space are scaled, and also both consider a sequence of initial
conditions $x_j$, such that after space-scaling, the ``relative initial
condition'' $x_j/a_j$ is a bounded distance away from the set $\frac
{n_j h}{a_j}\eq$. More precisely, we define the following property
which is common to both near and far fluid limit sequences. Thus by
assuming this property in the statement of Theorem~\ref
{thsubseqconv}, the theorem applies to both near and far fluid limit
sequences.

%The proof of Theorem~\ref{thsubseqconv} requires a few lemmas.
%Lemmas~\ref{lemFSLLN},~\ref{lemthin}, and~\ref{lemuniformint} all
%assume that we start with a sequence that satisfies the following
%properties:

%
%pr1 #&#
%
\begin{property}
\label{propa}
$\{ ({{\mathbf{x}}_{j}},a_j)\}$ is a sequence of initial condition $x_j$,
threshold factor $n_j$ and scale $a_j$ triples for which $a_j
\rightarrow\infty$. Moreover for each $j$, $n_j>0$, $a_j >0$ and some
closed, bounded $\eq\in\bar{\mathrm{X}}$,
%$ \frac{U^\ic{j}(0)}{a_j} \rightarrow\bar{U}(0)$, $\frac{V^
%and $\bar{V}(0) \geq0$. Moreover for each $j$,
%
\[
\frac{n_j}{a_j} \leq c_1 \quad\mbox{and}\quad \biggl\llVert
\frac{x_j}{a_j} \biggr\rrVert_{ {n_j h \eq}/{a_j} } \leq c_2 \qquad\mbox{for
some } c_1>0 \mbox{ and }c_2>0. %\mbox{ for some } c_1>0, \\
%some } c_2>0.\\%, \mbox{and} \\
% \frac{U^\ic{j}(0)}{a_j} &\rightarrow\bar{U}(0), \frac{V^
\]
%
%for some $\bar{U}(0)$
\end{property}
%
%This property ensures that there are well defined residual arrival and
%service times in the fluid limit.
%

%%%%%%%%%%%%%%%%%%%%%%%%%%%%%%%%%%%%%%
%
% THEOREM 1
%
%%%%%%%%%%%%%%%%%%%%%%%%%%%%%%%%%%%%%

%s3.1 #&#
\subsection{Convergence to a fluid limit along a subsequence}
\label{secfluidlimitconvergence}
The proof of the following theorem parallels the proof of Theorem 4.1
of Dai~\cite{Dai95}. However, the proof of our theorem differs in that
we require some specialized treatment for our fluid limit construction
and for the ingress discarding feature of the network. We state the
theorem here and present the proof in the \hyperref[app]{Appendix}.

%th1 #&#
%
\begin{theorem}
\label{thsubseqconv} Suppose $\{({{\mathbf{x}}_{j}},a_{j})\}$ is a
sequence satisfying
\label{referback}Property~\ref{propa} (on pa\-ge~\pageref{propa}).
%
% Suppose one of the following cases hold for some constant $\const$:
%
% \noindent\textbf{Near Fluid Limit (Near FL) Case:} $\icb{m}=
%satisfying:
% \begin{alignat*}{3}
%
% \n{\frac{x_m}{n_m}}{h \equilib} = \frac{\n{x_m}{n_m h \equilib}}{n_m}
%&\leq\const& \mbox{ and } && n_m \rightarrow\infty.
% %
% \intertext{\textbf{Far Fluid Limit (Far FL) Case:} $\icb{m}=
% condition and system scale pairs satisfying:}
%
% \n{ \frac{x_m}{n_m} }{h \equilib} = \frac{\n{x_m}{n_m h
% \equilib}\rightarrow\infty.
% \end{alignat*}
%
Then for almost all $\omega$ there exists a subsequence $\{({{\mathbf
{x}}_{m}},a_{m})\}
\subseteq\{({{\mathbf{x}}_{j}},a_{j})\}$ for which
\[
\frac{ \frak{X}^{{\mathbf{x}}_{m}}( a_mt)}{a_m} \rightarrow\bar{\frak{X}}(t)
\qquad\mbox{u.o.c.}
\]
for some fluid model trajectory $\bar{\frak{X}}(\cdot)$ with components
\[
\bar{\frak{X}}(\cdot) \triangleq\bigl[\bar{X}(\cdot);\bar{T}(\cdot
);\bar{
\Lambda}(\cdot);\bar{h}\bigr],
\]
where, in turn, the process $\bar{X}(\cdot)$ has components
\[
\bar{X}(\cdot) \triangleq\bigl[\bar{Q}(\cdot);\bar{U}(\cdot);\bar
{V}(\cdot);
\bar{H}(\cdot)\bigr],
\nonumber
\]
where $\bar{H}(\cdot)\equiv0$. The process $\bar{\frak{X}}(\cdot
)$ may
depend upon $\omega$ and the choice of subsequence $\{({{\mathbf
{x}}_{m}},a_{m})\}$ but
must satisfy the following properties for all $t\geq0$:
%
%e12 #&#
%e13 #&#
%e14 #&#
%e15 #&#
%e16 #&#
%e17 #&#
%e18 #&#
%e19 #&#
%
\begin{eqnarray}
\label{eqUV}
&\displaystyle \bar{U}_f(t) = \bigl(t-\bar{U}_f(0)
\bigr)^+,\qquad %\\ %\mbox{ For each $f\in
\bar{V}_k(t) = \bigl(t-
\bar{V}_k(0)\bigr)^+;&
\\
\label{eqTinc} &\displaystyle \bar{T}_k(t) \mbox{ is
nondecreasing and starts from zero};&
\\
\label{eqI} &\displaystyle \bar{I}_i(t):= t - C_i \bar{T}(t) \mbox{ is
nondecreasing};&
\\
\label{eqDeparturelimit1} &\displaystyle \bar{D}_k(t):= \mu_{s(k)}\bigl(
\bar{T}_k(t) - \bar{V}_k(0)\bigr)^+;&
\\
\label{eqArrivallimit} &\displaystyle \bar{A}(t):= P^\top\bar{D}(t) + \bar{
\Lambda}(t);&
\\
\label{eqQbar} &\displaystyle \bar{Q}(t):= \bar{Q}(0) + \bar{A}(t)-\bar{D}(t);&
\\
\label{eqQbig0limit} &\displaystyle \bar{Q}(t) \geq0;&
\\
\label{eqidlelimit} &\displaystyle \int_0^{\infty} \bigl(C
\bar{Q}(t)\bigr) \,d\bar{I}(t)=0,&
\end{eqnarray}
where (\ref{eqUV}), (\ref{eqTinc}) and (\ref{eqDeparturelimit1})
hold for each flow $f$ and class $k$, while (\ref{eqI}) holds for each
station $i$. Assignments (\ref{eqI}), (\ref{eqDeparturelimit1}),
(\ref
{eqArrivallimit}) and (\ref{eqQbar}) define $\bar{I}(t)$, $\bar
{D}(t)$, $\bar{A}(t)$ and $\bar{Q}(t)$, respectively. Also, the
following hold for each flow $f$ for regular $t\geq0$:
%
%&&\mbox{whenever $\bar{Q}_k(t)>\bar{h}$ for some $k\in\mathcal
%{C}(f)$},\\
%& &\mbox{whenever $\bar{Q}_k(t)<\bar{h}$ for all $k\in\mathcal
%{C}(f)$},\\
%%
%
%e20 #&#
%e21 #&#
%e22 #&#
%
\begin{eqnarray}
\label{eqHlimit} \bd{\Lambda}_f(t) &=& 0 \qquad\mbox{whenever $
\bar{Q}_k(t)>\bar{h}$ for some $k\in\mathcal{C}(f)$},
\\
\label{eqGlimit} \bd{\Lambda}_f(t) &=& \alpha_f 1
\bigl(t\geq\bar{U}_f(0)\bigr) \qquad\mbox{whenever $\bar{Q}_k(t)<
\bar{h}$ for all $k\in\mathcal{C}(f)$},
\\
\label{eqGlimit2} \bd{\Lambda}_f(t) &\leq&\alpha_f.
\end{eqnarray}
Also, for station $i$ and for any $k,l$ such that $ \{k,l\} \in C(i)$
the following properties are satisfied for all regular $t \geq0$:
%
%%
%&w_k^{-1}\dot{\bar{D}}_k(t) \geq w_l^{-1}\dot{\bar{D}}_l(t) \nonumber\\
%& \mbox{whenever $\bar{Q}_k(t)>0$ }, \\
%&w_k^{-1}\dot{\bar{D}}_k(t) = w_l^{-1}\dot{\bar{D}}_l(t) \nonumber\\
%&\mbox{whenever $\bar{Q}_k(t)>0$ and
%$Q_l(t)>0$}.
%%
%}
%
%e23 #&#
%e24 #&#
%
\begin{eqnarray}
\label{eqRRdrain2} w_k^{-1}\dot{\bar{D}}_k(t)
&\geq& w_l^{-1}\dot{\bar{D}}_l(t) \qquad
\mbox{whenever $Q_k(t)>0$},
\\
\label{eqRRdrain1} w_k^{-1}\dot{\bar{D}}_k(t)
&=& w_l^{-1}\dot{\bar{D}}_l(t) \qquad\mbox{whenever
$Q_k(t)>0$ and $Q_l(t)>0$}.
\end{eqnarray}
%
% In addition the following case-specific conditions on $\bar{h}$,
% and $\bar{X}(0)$ hold:
% %
% \begin{equation}
% \label{eqcase1Xhe}
% &\mbox{\textbf{Near FL Case:} }& \bar{h}&= h,&
% \n{\bar{X}(0)}{\bar{h}\equilib} &\leq\const;\\
% \label{eqcase2Xhe}
% &\mbox{\textbf{Far FL Case:} }& 0\leq\bar{h} &\leq h,&
% \n{\bar{X}(0)}{\bar{h}\equilib} &= 1.
% \end{equation}
\end{theorem}
See the \hyperref[app]{Appendix} for the proof. Next we state precisely the definitions
of a near fluid limit sequence and far fluid limit sequence that we
discussed earlier in Section~\ref{sproofoutline}. After defining these
sequences, we derive two corollaries to Theorem~\ref{thsubseqconv}
that apply to each of these types of sequences.
%
%de1 #&#
%
\begin{definition}[(Near fluid limit sequence)]
\label{defNFL}
$\{({{\mathbf{x}}_{j}},a_{j})\}$ is a \textit{near fluid limit
sequence} with respect
to a
closed, bounded $h \eq\in\bar{\mathsf{X}}$ if $a_j= n_j$,
$n_j\rightarrow\infty$ and
\[
\biggl\llVert\frac{x_j}{n_j} \biggr\rrVert_{ h \equilib} =
\frac{\llVert
x_j \rrVert_{ n_jh \equilib}}{n_j} \leq\zeta
\]
for each $j$ and for some $\zeta>0$.
\end{definition}
%
%de2 #&#
%
\begin{definition}[(Far fluid limit sequence)]
$\{({{\mathbf{x}}_{j}},a_{j})\}$ is a \textit{far fluid limit
sequence} with respect to a
closed, bounded $h \eq\in\bar{\mathsf{X}}$ if $a_j= n_j\llVert
\frac{x_j}{n_j} \rrVert_{ h \equilib}$, $a_j\rightarrow
\infty$ and
\[
\biggl\llVert\frac{x_j}{n_j} \biggr\rrVert_{ h \equilib} =
\frac{\llVert
x_j \rrVert_{ n_jh \equilib}}{n_j} > \zeta
\]
for each $j$ and for some $\zeta>0$.
\end{definition}

As was discussed earlier, the near fluid limit sequence is defined so
that the sequence of scaled initial conditions remains a bounded
distance away from the set $h\eq$ while the far fluid limit is defined
so that the sequence of scaled initial conditions is bounded away from
the set $h\eq$.

%%%%%%%%%
% COROLLARY 1
%%%%%%%%%
%
%co1 #&#
%
\begin{corollary}
Suppose that $\{({{\mathbf{x}}_{j}},a_{j})\}$ is a near fluid limit
sequence with respect
to a closed, bounded $h \eq\in\bar{\mathsf{X}}$.
% with $n_m \rightarrow\infty$, $a_m = n_m$, and
%for each $m$ and for some $\const>0$. %We call such a sequence a
%
Then for almost all $\omega$ there exists a subsequence $\{({{\mathbf
{x}}_{m}},a_{m})\}
\subseteq\{({{\mathbf{x}}_{j}},a_{j})\}$ for which
\[
\frac{ \frak{X}^{{\mathbf{x}}_{m}}( a_mt)}{a_m} \rightarrow\bar{\frak
{X}}(t)\qquad\mbox{u.o.c.},
\]
where $\frak{X}(\cdot)$ satisfies fluid model equations (\ref
{eqUV})--(\ref{eqRRdrain1}). Moreover
\[
\bar{h} = h \quad\mbox{and}\quad \bigl\llVert\bar{X}(0) \bigr\rrVert_{ \bar
{h}\eq}
\leq\zeta.
\]
\end{corollary}
\begin{pf}
The discarding thresholds before scaling are $n_jh$, and thus
after scaling they are $n_jh/ a_j= h$ for each $j$. Thus
$\bar{h}=h$. Also $a_j\rightarrow\infty$ and \mbox{$ \llVert x_j /a_j
\rrVert_{ h \eq} \leq\zeta$}, and\vadjust{\goodbreak} thus the sequence $\{({{\mathbf
{x}}_{j}},a_{j})\}$
satisfies Property~\ref{propa}. By Theorem~\ref{thsubseqconv} there
exists a subsequence $\{({{\mathbf{x}}_{m}},a_{m})\}$ such that
${\frak{X}^{{\mathbf
{x}}_{m}}(
a_mt) }/{a_m}$ converges u.o.c. to a fluid trajectory
satisfying (\ref{eqUV})--(\ref{eqRRdrain1}). By Theorem \ref
{thsubseqconv}, the subsequence $x_m/a_m$ converges to an
initial state of the fluid trajectory $\bar{X}(0)$. Since $\llVert
x_m /a_m \rrVert_{ h \equilib} \leq\zeta$, it must be that $\llVert
\bar{X}(0) \rrVert_{ \bar{h} \eq} \leq\zeta$.
\end{pf}

%%%%%%%%%
% COROLLARY 2
%%%%%%%%%
%
%co2 #&#
%
\begin{corollary}\label{cor2}
Suppose that $\{({{\mathbf{x}}_{j}},a_{j})\}$ is a far fluid limit
sequence with respect
to a closed, bounded $h \eq\in\bar{\mathsf{X}}$. %with $a_m = n_m
% \[ \n{ \frac{x_m}{n_m} }{h \equilib} = \frac{\n{x_m}{n_m h
% for each $m$ and for some $\const>0$.
%
Then for almost all $\omega$ there exists a subsequence $\{({{\mathbf
{x}}_{m}},a_{m})\}
\subseteq\{({{\mathbf{x}}_{j}},a_{j})\}$ for which
\[
\frac{ \frak{X}^{{\mathbf{x}}_{m}}( a_mt)}{a_m} \rightarrow\bar{\frak
{X}}(t)\qquad
\mbox{u.o.c.},
\]
where $\frak{X}(\cdot)$ satisfies fluid model equations (\ref
{eqUV})--(\ref{eqRRdrain1}). Moreover
\[
\bar{h} \in[0,h/\zeta] \quad\mbox{and}\quad \bigl\llVert\bar{X}(0) \bigr
\rrVert_{ \bar{h}\eq} = 1.
\]
\end{corollary}
\begin{pf}
Note that
\[
\biggl\llVert\frac{x_j} {a_j} \biggr\rrVert_{ {n_jh\eq}/{a_j} } =
\frac
{n_j}{a_j} \biggl\llVert\frac{x_j} {n_j} \biggr\rrVert_{ h \eq} =
1 % &= \n{ \frac{x_m} {n_m \n{ \frac{x_m}{n_m} }{h \equilib} }}{
%h}{a_m} \eq} \\
% &= \n{ \frac{ \frac{x_m}{ n_m} - h e_m + h e_m} { \n{ \frac{x_m}{n_m}
%}{h \equilib} }}{\frac{n_m h}{a_m} \eq}\\
% &= \n{ v_m + \frac{ h}{ \n{ \frac{x_m}{n_m} }{h \equilib} } e_m }{
% &= \n{ v_m + \frac{n_m h }{ n_m \n{ \frac{x_m}{n_m} }{h \equilib} }
%e_m}{\frac{n_m h}{a_m} \eq}\\
% &= \n{ v_m + \frac{n_m h }{ a_m } e_m}{\frac{n_m h}{a_m} \eq}\\
% &= 1
\]
for each $j$.
This combined with the fact that $a_j\rightarrow\infty$ implies
that $\{({{\mathbf{x}}_{j}},a_{j})\}$ satisfies Property \ref
{propa}. By Theorem \ref
{thsubseqconv} there exists a subsequence $\{({{\mathbf
{x}}_{m}},a_{m})\}$ such that
${\frak{X}^{{\mathbf{x}}_{m}}( a_mt) }/{a_m}$
converges u.o.c. to a
fluid trajectory satisfying (\ref{eqUV})--(\ref{eqRRdrain1}). The
above equation also implies that $\llVert \bar{X}(0) \rrVert_{ \bar
{h} \eq}=1$. The
subsequence of scaled thresholds satisfies $n_mh/ a_m= h / \llVert
{x_m}/{n_m} \rrVert_{ h \eq} < h/\zeta$. By Theorem~\ref{thsubseqconv} the
subsequence $n_mh/a_m$ converges, and the convergence must be
to a number in the range $[0, h/\zeta]$ because of the preceding
inequality relation.
\end{pf}

%s3.2 #&#
\subsection{Convergence along subsequences to
convergence along sequences}
%Theorem~\ref{thmfluidconv} shows that the fluid limit, either
In the previous section, we showed that for both near and far fluid
limit sequences, we can extract a sample path dependent subsequence
that converges to a fluid model trajectory. The objective of this
section is to use this subsequence result to show convergence of a
functional of the original sequence. In particular, we show in
Lemma~\ref{lemsubseqtoseq} that if a functional $\frak{F}$ of any
fluid model trajectory goes to zero in a time not more than a constant
times the scaled initial condition's distance from $\bar{h}\eq$, then
the value of that functional applied to the fluid limit sequence of
trajectories converges almost surely. In later sections, we will invoke
Lemma~\ref{lemsubseqtoseq} choosing $\frak{F}$ to extract the
service rates from the fluid model, and later choosing $\frak{F}$ to
extract the distance from a set $h \eq$. Lemma~\ref{lemsubseqtoseq}
is a generalization of an argument used by Dai in
the proof of Theorem 4.2 of~\cite{Dai95}.
%%%%%%%%%%%%%%%%%%%%%%%%%%%%%%%%%%%%%%%%%%%%%%
% Lemma 4
% 4444444444444444444444444444444444444444444
%
%%%%%%%%%%%%%%%%%%%%%%%%%%%%%%%%%%%%%%%%%%%%%%%%%
%
%le1 #&#
%
\begin{lemma}
\label{lemsubseqtoseq}
Suppose that $\frak{F}$ is a functional that maps $\R^{r} \times\R^+$
into $\R^{s} \times\R^+ $ where $r$ is the dimension of $\frak
{X}^{\mathbf{x}}(\cdot)$ and $s$ is arbitrary. Also suppose
that $\frak{F}$ is
continuous on the topology of uniform convergence on compact sets. If
the following is true:
\begin{itemize}
\item The fluid model equations (\ref{eqUV})--(\ref{eqRRdrain1}) are
such that for any trajectory $\bar{\frak{X}}(\cdot)$ and $\bar{h}
\geq
0$ that satisfies them, there exists some closed bounded $\eq\in
\bar{\mathsf
{X}}$ for which
%
%e25 #&#
%
\begin{equation}
\label{eqcase1F} \frak{F}\circ\bigl[\bar{\frak{X}}(\cdot)\bigr](t)
\equiv0\qquad
\forall t\geq t_0 \bigl\llVert\bar{X}(0) \bigr
\rrVert_{ \bar
{h}\equilib}.
\end{equation}
\end{itemize}
Then, for any sequence $\{({{\mathbf{x}}_{j}},a_{j})\}$ satisfying
Property~\ref{propa}
where the relation \label{corr2}$\llVert x_j/a_j \rrVert_{
{n_j h\eq}/{a_j} } \leq c$ of Property~\ref{propa} is satisfied with
constant $c>0$,
%
%e26 #&#
%
\begin{equation}
\label{eqseqconclude} \biggl\llvert\frak{F}\circ\biggl[ \frac{1}{\s{j}}
\frak{X}^{{\mathbf
{x}}_{j}}(\s{j}\cdot) \biggr](t) \biggr\rrvert\rightarrow0
\qquad\mbox{a.s.}
\end{equation}
for each $t \geq c t_0$.
\end{lemma}
\begin{pf}
By Theorem~\ref{thsubseqconv}, for almost all sample paths $\omega$,
and for any subsequence $\{({{\mathbf{x}}_{m}},a_{m})\}\subseteq\{
({{\mathbf{x}}_{j}},a_{j})\}$ there is a
sample-path-dependent further-subse\-quence $\{({{\mathbf{x}}_{r(\omega
)}},a_{r(\omega)})\}
\subseteq\{({{\mathbf{x}}_{m}},a_{m})\}$ for which
\[
\frac{\frak{X}^{{\mathbf{x}}_{r(\omega)}}(\s{r(\omega
)}t,\omega)
}{\s{r
(\omega)}}\rightarrow\bar{\frak{X}}(t,\omega)\qquad
\mbox{u.o.c.},
\]
where $\bar{\frak{X}}(t,\omega)$ satisfies (\ref{eqUV})--(\ref
{eqRRdrain1}) as well as \label{correction}$\llVert \bar{X}(0)
\rrVert_{ \bar{h}\equilib}\leq c$ since each $x_j/a_j$ has a
distance from
$\frac{n_jh}{a_j} \eq$ that is no more than $c$ by the lemma's
assumption. The notation $r(\omega)$ and $\bar{\frak
{X}}(t,\omega)$
emphasize that the further-subsequence and fluid trajectory depend on
$\omega$. Now fix an $\omega$ for which subsequences have convergent
further subsequences as described. For the next few steps we suppress
the $\omega$ arguments to simplify notation.
% \begin{equation}
% \label{eqfrak}
% \frak{F}\circ\left[\bar{\frak{X}}(\cdot)\right](t) = 0
% \mbox{for each $t\geq\n{ \bar{X}(0)}{\bar{h} \equilib} \cdot
%t_0$}. % |X(0) - h\ul{e}|t_0$}.
% \end{equation}
Because $\frak{F}$ is assumed to be continuous on the topology of
uniform convergence on compact sets, we have %from (\ref{eqcase1F})
\[
\frak{F} \circ\biggl[ \mbox{$ \frac{ \frak{X}^{{\mathbf
{x}}_{r}}(\s
{r}
\cdot)}{\s{r}}$} \biggr](t) \rightarrow
\frak{F} \circ\bigl[\bar{\frak{X}}(\cdot)\bigr](t) \qquad\mbox{u.o.c.}
\]
Consequently,
\[
% \label{eqeek}
\biggl\llvert\frak{F} \circ\biggl[\mbox{$ \frac{ \frak{X}^{{\mathbf
{x}}_{r}}(\s{r
} \cdot)}{\s{r}}$}
\biggr](t) \biggr\rrvert\rightarrow0 %
\]
for each $t\geq c t_0$. So for this fixed $\omega$, any subsequence
$\{({{\mathbf{x}}_{m}},a_{m})\}\subseteq\{({{\mathbf
{x}}_{j}},a_{j})\}$ has a further subsequence $\{({{\mathbf{x}}_{r
(\omega)}},a_{r (\omega)})\}\subseteq\{({{\mathbf{x}}_{m}},a_{m})\}
$ for which the above holds. Therefore the
original sequence $\{({{\mathbf{x}}_{j}},a_{j})\}$ converges for this
fixed $\omega$. The
same argument can be used to conclude that this holds for almost all~$\omega$.
Thus, we have~(\ref{eqseqconclude}).
\end{pf}
\subsection{Convergence to fluid model rates on a compact time interval}
\label{sratescompact}

The objective of this section is to use
Lemma~\ref{lemsubseqtoseq} to conclude that the rates of the
stochastic system are close to those of the fluid model over a finite
time interval. It will remain to show that the rates are close over the
long-term.\vadjust{\goodbreak}
%%%%%%%%%%%%%%%%%%%%%%%%%%%%%%%%%%%%%%%%%%%%%%
% THEOREM 2
% 2222222222222222222222222222222222
%
%%%%%%%%%%%%%%%%%%%%%%%%%%%%%%%%%%%%%%%%%%%%%%%%%
%
%th2 #&#
%
\begin{theorem}
\label{thratecompactset}
Suppose there exists $t_0>0$, a closed, bounded $\eq\in\bar
{\mathsf{X}}$, and rate vector $R\in\mathbb{R}^K_+$ such that
%
%e27 #&#
%
\begin{equation}
\label{eqth30} M^{-1}\bd{T}(t) \equiv R\qquad\forall t\geq
t_0\bigl\llVert\bar{X}(0) \bigr\rrVert_{ \bar{h}\eq}
\end{equation}
for any fluid model trajectory $\bar{\frak{X}}(\cdot)$ and $\bar{h}>0$
that satisfies (\ref{eqUV})--(\ref{eqRRdrain1}). Then for any
positive $\gamma<1$ and $\zeta<1$ there exists $L_1(\zeta,\gamma)$ such
that for all $n\geq L_1$,
%
%e28 #&#
%
\begin{equation}
\label{eqrc4} \inf_{\llVert {x}/{n} \rrVert_{ h\equilib}\leq\zeta
} \E\bigl[ M^{-1} T^{{\mathbf{x}}_{}}(
n t_0) \bigr] \geq R (1-\zeta) (1-\gamma)n t_0.
\end{equation}
\end{theorem}
%
%%%%%%%%%%%%%%%%%%%%%%%%%%%%%%%%%%%%%%%%%%%%%%
% THEOREM 2 PROOF
% 2222222222222222222222222222222222
%
%%%%%%%%%%%%%%%%%%%%%%%%%%%%%%%%%%%%%%%%%%%%%%%%%
%
\begin{pf}
Let $\{({{\mathbf{x}}_{j}},a_{j})\}$ be a near fluid limit sequence: a
sequence of
threshold scale and initial condition pairs satisfying $a_j
=n_j
\rightarrow\infty$ and $\llVert x_j/\break n_j \rrVert_{ h\equilib} \leq
\zeta
$. We
invoke Lemma~\ref{lemsubseqtoseq} by picking $\frak{F}$ so that
\[
\frak{F} \circ\bigl[\bar{\frak{X}}( \cdot)\bigr](t) %\bar{F}(t) &=
:=
\bar{T}\bigl(\zeta^{-1} t\bigr)-\bar{T}(t) - MR\bigl(
\zeta^{-1}-1\bigr)t.
\]
$\frak{F}$ is easily seen to be continuous on the topology of uniform
convergence on compact sets. Also note that $\frak{F} \circ[\bar
{\frak
{X}}( \cdot)](t)=0$ for all $t\geq t_0\llVert \bar{X}(0) \rrVert_{
\bar{h}\eq}$ by
(\ref{eqth30}). By Lemma~\ref{lemsubseqtoseq},
\[
\lim_{j\rightarrow\infty} \biggl\llvert\frac{ T^{{\mathbf
{x}}_{j}}(n_jt_0 )-T^{{\mathbf
{x}}_{j}}( \zeta
n_j
t_0)}{n_j(1-\zeta)t_0 } -MR \biggr
\rrvert= 0 \qquad\mbox{a.s.},
\]
where we have used the fact that $\llVert x_j/n_j \rrVert_{
h\equilib}
\leq
\zeta$ to choose the $c$ of Lemma~\ref{lemsubseqtoseq} to be
$\zeta$
and selected $t=\zeta t_0$. The left-hand side of the above identity
is bounded from above by a constant for all $j$, and thus by the
dominated convergence theorem~\cite{Durrett},
%
%e29 #&#
%
\begin{equation}
\label{eqcat3} \lim_{j\rightarrow\infty} \E\biggl\llvert\frac{
T^{{\mathbf{x}}_{j}}(n_jt_0 )-T^{{\mathbf
{x}}_{j}}( \zeta
n_j
t_0)}{n_j(1-\zeta)t_0 } -MR
\biggr\rrvert= 0.
\end{equation}
Also note (\ref{eqcat3}) holds for any sequence $\{({{\mathbf
{x}}_{j}},a_{j})\}$ with
$n_j= a_j\rightarrow\infty$ and $\llVert x_j/n_j \rrVert_{
h\equilib}\leq\zeta$, because these were the only restrictions for
our initial choice of sequence.

Now pick a positive constant $\gamma<1$. Observe that there exists a
constant $L_1(\gamma,\zeta)$ such that whenever $n> L_1$,
\[
\inf_{\llVert x/n \rrVert_{ h\equilib} \leq\zeta} \frac{\E[
T^{\mathbf{x}}( n t_0) - T^{\mathbf{x}}( n \zeta
t_0) ]}{n
(1-\zeta)t_0} \geq MR(1-\gamma)
\]
for if otherwise we could construct a sequence $\{({{\mathbf
{x}}_{j}},a_{j})\}$ that
violates (\ref{eqcat3}). By the monotonicity of $T^{{\mathbf
{x}}_{j}}
(\cdot
)$, we have (\ref{eqrc4}).
\end{pf}

%s3.4 #&#
\subsection{Stochastic system attracted to $h\eq$}
\label{secFluidLimittoEquil}

The objective of this section is to show that the scaled state of the
stochastic system is attracted to $h\eq$. In particular we show that
the scaled state's expected distance from $h\eq$ declines
geometrically
(roughly) for starting scaled states outside a neighborhood of $h\eq$.
Since the proof technique is similar that of Theorem 3.1 of Dai \cite
{Dai95} we choose to provide the proof in the \hyperref[app]{Appendix}.\vadjust{\goodbreak}

%Before stating the theorem of this section, we define:
% \label{eqYdef}
% y:= x/n, \y=(y,n), Y^{\y}(t):= \frac{1}{n}X^\ic{}(nt).

%%%%%%%%%%%%%%%%%%%%%%%%%%%%%%%%
%
% Theorem 3
% 333333333333333333333333333333
%%%%%%%%%%%%%%%%%%%%%%%%%%%%%%%%%%%
%
%th3 #&#
%
\begin{theorem}
\label{thPushtoEquil}
Suppose that there exists $t_0>0$ and a closed, bounded $\eq\in
\bar{\mathsf
{X}}$ such that
%
%e30 #&#
%
\begin{equation}
\label{eqth20} \bigl\llVert\bar{X}(t) \bigr\rrVert_{ \bar{h}\equilib}
\equiv0\qquad
\forall t\geq t_0\bigl\llVert\bar{X}(0) \bigr
\rrVert_{ \bar{h}\equilib}
\end{equation}
for any fluid model trajectory $\bar{\frak{X}}(\cdot)$ and $\bar
{h}\geq
0$ that satisfies (\ref{eqUV})--(\ref{eqRRdrain1}). Then the
following conclusions are true:
\begin{longlist}
\item For any $\zeta>0$, and any positive $\delta<1$ there exists
$L_2(\zeta,\delta)$ such that for all $n \geq\zeta^{-1} L_2$ and all
$x\dvtx  \llVert x/n \rrVert_{ h\eq} > \zeta$,
\[
\E\biggl\llVert\frac{1}{n} X^{{\mathbf{x}}_{}} \biggl( n t_0
\biggl\llVert\frac{x}{n} \biggr\rrVert_{ h \eq} \biggr) \biggr
\rrVert_{ h \eq} \leq\delta\biggl\llVert\frac{x}{n} \biggr
\rrVert_{ h \eq}. %\E\n{ Y^{\y}( t_0 \n{y}{h\eq})}{h\eq} \leq\delta
\]
\item For any $\zeta>0$, and any $b>0$ there exists $L_3(\zeta,b)$
such that for all $n\geq L_3$ and all $x\dvtx  \llVert x/n \rrVert_{ h\eq
} \leq\zeta$,
\[
\E\biggl\llVert\frac{1}{n} X^{{\mathbf{x}}_{}} (n t_0) \biggr
\rrVert_{ h \eq
} \leq b. % \E\n{ Y^{\y}( t_0 )}{h\eq} \leq b
\]
\end{longlist}
\end{theorem}
See the \hyperref[app]{Appendix} for the proof.

The objective of the next lemma is to show that the results of
Theorem~\ref{thPushtoEquil} imply that the expected return time of the
scaled state to the $\zeta$ ball around $h\eq$ is small. The proof of
Lemma~\ref{lemcontract} is adapted from the proof of Theorem~2.1(ii)
of~\cite{MeynTweedieStateDep}, which was for a discrete time Markov
chain. Since the lemma is an adaptation of a previous result, we \label
{provide}provide the proof in the \hyperref[app]{Appendix}.

%At this point we invoke Theorem 3.1 of~\cite{Dai95}, which requires
%that the interarrival times are unbounded and spread out, as stated in
%(\ref{equnbounded}) and (\ref{eqspreadout}). Since the lemma is of
%the similar to Theorem 2.1(ii) of~\cite{MeynTweedieStateDep}, we
%choose to state the Lemma here and provide the proof in the appendix.

%In Lemma~\ref{lemcontract} and subsequently, to express $Y^{\y}(t)$
%without specifying an initial condition we write $Y^n(t)$ where the
%choice of initial condition is implicit by the choice of probability
%measure. We define $\Pb_y$ to be a probability measure for which
% $ \Pb_y\{Y^n(0) = y \}=1, $
%and thus,
% $ Y^n(t) = Y^{\y}(t) \mbox{$P_y$-a.s.} $
%and
% $ \E_y[Y^n(t)] = \E[Y^{\y}(t)].$

%%%%%%%%%%%%%%%%%%%%%%%%%%%%
%
% LEMMA STOPPING TIME
%
%%%%%%%%%%%%%%%%%%%%%%%%%%%%%
%
%le2 #&#
%
\begin{lemma}
\label{lemcontract}
Suppose (\ref{equnbounded}) and (\ref{eqspreadout}) are satisfied
and for some $n>0$, $h\geq0$, and a closed, bounded $\eq\in\bar{\mathsf
{X}}$ we have % the following inequalities:
%
%e31 #&#
%e32 #&#
%
\begin{eqnarray}
\label{eqlyaplem} \E\biggl\llVert\frac{1}{n} X^{{\mathbf{x}}_{}} \biggl
( n
t_0 \biggl\llVert\frac{x}{n} \biggr\rrVert_{ h \eq}
\biggr) \biggr\rrVert_{ h \eq} &\leq& \delta\biggl\llVert\frac{x}{n}
\biggr\rrVert_{ h \eq}\qquad \forall x\dvtx  \llVert x/n \rrVert_{ h\eq} >
\zeta,
\\
%&\E_y\n{Y^{n} (t_0 |y |) }{h\eq} &&\leq\delta\n{y}{h\eq} &
\label{eqlyaplem2} \E\biggl\llVert\frac{1}{n}
X^{{\mathbf{x}}_{}} (n t_0) \biggr\rrVert_{ h \eq
} &\leq& b\qquad
\forall x\dvtx  \llVert x/n \rrVert_{ h\eq} \leq\zeta. %&\E_y\n{Y^{n} (t_0 )
%} {h\eq} &&\leq b &\mbox{for all $\n{y}{h
\end{eqnarray}
Then $X$ is positive Harris recurrent and
%
%e33 #&#
%
\begin{equation}
\label{eqtau1lem} \sup_{x \in B} \E_x\bigl[
\tau^n_B( n t_0)\bigr] \leq n
t_0 \biggl[1 + \frac
{\zeta
+b}{1-\delta} \biggr],
\end{equation}
where $B \triangleq\{ x\dvtx  \llVert x/n \rrVert_{ h\eq} \leq\zeta\}
$ and $\tau^n_B(n
t_0)$ is
defined by
%
%e34 #&#
%
\begin{equation}
\label{eqtaulem} \tau^n_B(n t_0)
\triangleq\inf\bigl\{t\geq n t_0\dvtx  X^{n}( t) \in B \bigr\}.
\end{equation}
\end{lemma}
See the \hyperref[app]{Appendix} for the proof.

%s3.5 #&#
\subsection{Convergence of long-term rates}
\label{slongtermrates}

\label{sketch} 
\begin{figure}

\includegraphics{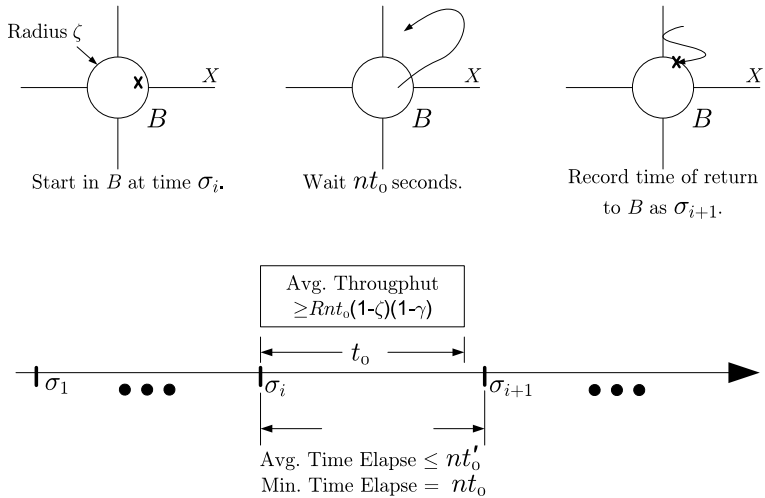}

\caption{The top half of the figure illustrates the definition of the stopping
times $\sigma_{i}, \sigma_{i+1},\ldots\,$. The bottom half illustrates the
intuition behind the proof of Theorem \protect\ref{thbigtheorem} by plotting
the stopping times on a time line, and showing the bound on expected
throughput between such stopping times.} %That the expected time elapse
%between $\sigma_i$ and $\sigma_{i+1}$ is upper-bounded by $t'_0$ is a
%consequence of Lemma~\ref{lemcontract}. The lower-bound on expected
%throughput between stopping times comes from Theorem
\label{figball}
\end{figure}

The objective of this section is to tie together all of
the preceding results to conclude in Theorem~\ref{thbigtheorem} that
the long-term rates of the stochastic\vadjust{\goodbreak} system are close to the fluid
rates for large enough $n$. First we pick $n$ large enough so that the
conclusions of Theorems~\ref{thratecompactset},~\ref{thPushtoEquil}
and Lemma~\ref{lemcontract} apply. Theorem~\ref{thratecompactset}
says that the stochastic system's rates are close to the fluid rates
for the first $nt_0$ seconds after having started with a scaled initial
condition $x/n$ in a $\zeta$-neighborhood of $h \eq$. To make a
conclusion about the long-term, we need to show that stochastic system
spends relatively little time away from the neighborhood in which
Theorem~\ref{thratecompactset} applies. Lemma~\ref{lemcontract}
tells us that the expected first return time of $X/n$ to a $\zeta
$-neighborhood of $h \eq$ that happens after $nt_0$ seconds is no more
than a constant times $nt_0$. Moreover, this constant can be made
arbitrarily small by picking $n$ larger. This argument is illustrated
by Figure~\ref{figball}. To formalize the argument we construct a
sequence of stopping times that occur on the first visit of $X/n$ to
%
%f4 #&#
%
the $\zeta$-neighborhood of $h \eq$ that occurs at least $nt_0$ seconds
after the last stopping time. We define random vectors $\rho_i$ that
track the cumulative service, divided by average service times, between
stopping times and relate these to the desired rate vector $R$ using
Theorem~\ref{thratecompactset}. We use ergodicity to argue that the
long-term average rates exist, and that this long-term limit must equal
the product of the expected value of $\rho_i$ times the lim inf of
$t/N(t)$ the inverse of the arrival rate of stopping times. Due to
Lemma~\ref{lemcontract}, this later quantity has an upper bound of
$nt_0$ times a constant that can be made small.

%%%%%%%%%%%%%%%%%%%%%%%%%%%%%%%%%%%%%%%%%%%%%%%%
%
% THEOREM 4
% 44444444444444444444444444444444444
%
%%%%%%%%%%%%%%%%%%%%%%%%%%%%%%%%%%%%%%%%%%%%%%%%%
%
%th4 #&#
%
\begin{theorem}
\label{thbigtheorem}
Suppose for some $t_0 > 0$ and some closed, bounded $\eq\in\bar{\mathsf
{X}}$ both of the following are true:
\begin{itemize}
\item For any fluid model trajectory $\bar{\frak{X}}(\cdot)$ and
$\bar
{h} \geq0$ that satisfies (\ref{eqUV})--(\ref{eqRRdrain1}),
%
%e35 #&#
%
\begin{equation}
\label{eqth31} \bigl\llVert\bar{X}(t) \bigr\rrVert_{ \bar{h} \eq}
\equiv0\qquad
\forall t\geq t_0\bigl\llVert\bar{X}(0) \bigr
\rrVert_{ \bar{h} \eq}.
\end{equation}
\item For any fluid model trajectory $\bar{\frak{X}}(\cdot)$ and
$\bar
{h} > 0$ that satisfies (\ref{eqUV})--(\ref{eqRRdrain1}),
%
%e36 #&#
%
\begin{equation}
\label{eqth32} M^{-1} \bd{T}(t) \equiv R\qquad \forall t\geq
t_0\bigl\llVert\bar{X}(0) \bigr\rrVert_{ \bar
{h}\eq},
\end{equation}
where $R\in\mathbb{R}^K_+$.
\end{itemize}
Then for any $\epsilon>0$, there exists a $n_c>0$ such that for all
$n\geq n_c$,
\[
\lim_{t \rightarrow\infty} \frac{D^{\mathbf{x}}(t)}{t} \geq(1-\epsilon)R
\qquad\mbox{a.s.}
\]
\end{theorem}
%
%%%%%%%%%%%%%%%%%%%%%%%%%%%%%%%%%%%%%%%%%%%%%%%%
%
% THEOREM 4 PROOF
% 44444444444444444444444444444444444
%
%%%%%%%%%%%%%%%%%%%%%%%%%%%%%%%%%%%%%%%%%%%%%%%%%
%
\begin{pf}
We observe that equations (\ref{eqth32}) and (\ref{eqth31}) are the
necessary conditions to apply Theorems~\ref{thratecompactset} and
\ref{thPushtoEquil}, respectively. Therefore, we may arbitrarily pick
the constants $\zeta$, $\delta$ and $b$ of Theorem~\ref{thPushtoEquil}
and the constants $\zeta$ and $\gamma$ of Theorem \ref
{thratecompactset} (using the same $\zeta$ value in Theorems \ref
{thratecompactset} as we use when we apply Theorem \ref
{thPushtoEquil}), and then fix an $n$ satisfying
%
%e37 #&#
%
\begin{equation}
\label{eqnchoice} n> \max\bigl[L_1(\zeta,\gamma),
\zeta^{-1}L_2(\zeta,\delta), L_3(\zeta,b)
\bigr]
\end{equation}
so that the conclusions of both Theorems~\ref{thPushtoEquil} and \ref
{thratecompactset} hold.

In addition, conclusions (i) and (ii) of Theorem~\ref{thPushtoEquil}
allow us to invoke Lemma~\ref{lemcontract} to complete (\ref
{eqtau1lem}) where $\tau_B^n( n t_0)$ is defined by (\ref
{eqtaulem}). Because the constants $\zeta$, $b$, $\delta$ can be
chosen arbitrarily, equations (\ref{eqtau1lem}) and (\ref
{eqnchoice}) imply that the ratio of the expected first hitting time
of $B$ ($n t_0$ seconds after having started in $B$) to $n t_0$ can be
made to be close to $1$ by choosing $n$ large enough. We collect some
of the constants in (\ref{eqtau1lem}) in the term $t_0'$ defined by
%
%e38 #&#
%
\begin{equation}
\label{eqvarphi} n t_0' = n t_0 \biggl[
1 +\frac{\zeta+b}{1-\delta} \biggr].
\end{equation}

We have also chosen $n$ large enough so that the following conclusion
from Theorem~\ref{thratecompactset} holds:
%
%e39 #&#
%
\begin{equation}
\label{eqRateconv} \inf_{\llVert x/n \rrVert_{ h\eq} \leq\zeta
} \E\bigl[ T^{\mathbf{x}}( n
t_0) \bigr] \geq M R (1-\zeta) (1-\gamma)n t_0.
\end{equation}
%
%The notation $T^{\y}( n t_0)$ above signifies $T^{\ic{}}(n t_0)$ where
%$\ic{}=(ny,n)$.
%%%%%%%%%%%%%%%%%%%%%%%%%%%%%%%%%%%%%%%%%%%%%%%%
%
% THEOREM 4 PROOF --- STOPPING TIME ARGUMENT
% 44444444444444444444444444444444444
%
%%%%%%%%%%%%%%%%%%%%%%%%%%%%%%%%%%%%%%%%%%%%%%%%%
%
Define the stopping times
%
%e40 #&#
%
\begin{equation}
\label{eqstoptimes} \sigma_0 = 0,\qquad %\nonumber\\
\sigma_{i+1} = \inf\bigl\{t \geq n t_0+
\sigma_{i}\dvtx X(t)\in B \bigr\}\qquad \forall i\geq0.
\end{equation}
Figure~\ref{figball} illustrates how these stopping times are defined.
Note that for any initial condition $x \in\mathrm{X}$ (the state space
of $X^n$) and index $i \geq1$,
%
%e41 #&#
%
\begin{equation}
\label{eqsigmabound1} \E_x[ \sigma_{i+1} -
\sigma_{i}] % = \E_{Y^{\y}(\sigma_{i})}[\tau_B^{n}(t_0) ]
\leq\sup_{\tilde{x}
\in B}\E_{\tilde{x}}
\bigl[\tau_B^n(nt_0)\bigr] \leq n
t_0'.
\end{equation}
This follows from the fact that $X^{\mathbf{x}}(\sigma_i)\in
B$, the strong
Markov property, the stopping time definitions (\ref{eqtaulem}) and
(\ref{eqstoptimes}) and expressions (\ref{eqtau1lem}) and (\ref
{eqvarphi}). Also, $X$ is positive Harris recurrent by Lemma \ref
{lemcontract} and therefore,
$\E_x[\sigma_1]< \infty$
for any $x \in\mathrm{X}$. We define a counting process $N(t)$ for the
stopping times $\sigma_i$ as
$
N(t) = \inf\{i\dvtx \sigma_i \leq t\}.
$
Because $X$ is positive Harris recurrent, $\sigma_i<\infty$ almost
surely, and therefore
$N(t) \rightarrow\infty$ a.s.
We now turn to bounding the expected ``arrival'' rate of the stopping
times $\sigma_i$. By (\ref{eqsigmabound1}) for each $i$,
%
%e42 #&#
%
\begin{equation}
\label{eqsigmabound3} \frac{\E_x[ \sigma_i ]}{ i} = \frac{ \sum
_{j=1}^{i -1} \E_x[\sigma_{j+1} - \sigma_j] + \E_x\sigma_1}{i} \leq n
t_0'(1- 1/i) + \frac{\E_x\sigma_1}{i}.
\end{equation}
Additionally, along any sample path
\[
\frac{t}{N(t)} \leq\frac{\sigma_{N(t)+1}}{N(t)+1} \frac{N(t)+1}{N(t)}.
\]
Thus by taking $\liminf_{t\rightarrow\infty}E_x(\cdot)$ of both sides,
and using (\ref{eqsigmabound3}) we have
\[
\liminf_{t \rightarrow\infty} \E_x \biggl[ \frac{t}{N(t)} \biggr]
\leq
n t_0'.
\]
Moreover, by Fatou's lemma
%
%e43 #&#
%
\begin{equation}
\label{eqfatou} \E_x \biggl[ \liminf_{t \rightarrow\infty}
\frac{t}{N(t)} \biggr] \leq\liminf_{t \rightarrow\infty} \E_x \biggl[
\frac{t}{N(t)} \biggr] \leq n t_0'.
\end{equation}

%Recall that the process $T^\ic{}(t)=T^{\y}(t)$ is defined in terms of
%the time scale of $X^n$, and not that of $Y^n$, which we defined in
%expression (\ref{eqYdef}) by scaling time by a factor of $n$.
%Therefore we define a re-scaled service process $\widetilde{T}{}^{\y}$
%with time re-scaled to match the time scale of the process $Y(t)$
%according to the definition
% $ \widetilde{T}{}^{\y}(a,b) = T^{\y}(nb) - T^{\y}(na). $
%We also
%
We define the random vectors $\rho_i= M^{-1} [{T}^{n}(\sigma_i +
\sigma_{i+1} ) - {T}^{n}(\sigma_i) ]$ to track the service between stopping
times $\sigma_i$.
%The $y$ superscript in $\widetilde{T}{}^{\y}$ is omitted because the
%initial condition will be specified implicitly by the choice of
%probability measure.
Note that for $i\geq1$ and each $x\in\mathrm{X}$,
%
%e44 #&#
%
\begin{equation}
\label{eqErho} \E_x[ \rho_i ] %\geq\E_y[ \tilde{\rho}_i ]
%&\geq\E_{Y^{\y}(\sigma_i)}[ M^{-1}\widetilde{T}{}^{n}( t_0 ) ] %
\geq\inf_{\tilde{x} \in B}\E_{\tilde{x}}\bigl[M^{-1}{T}^{\tilde
{x},n}(nt_0)
\bigr] %\nonumber\\
\geq R n t_0(1-\zeta) (1-\gamma).
\end{equation}
This follows from the fact that $X^{\mathbf{x}}(\sigma_i)\in
B$, the strong
Markov property, the definition of $\sigma_i$ (\ref{eqstoptimes}), the
definition of $\rho_i$, and relation (\ref{eqRateconv}). Figure \ref
{figball} illustrates the fact that the throughput between stopping
times $\sigma_i$ and $\sigma_{i+1}$ is lower-bounded according to
relation (\ref{eqErho}).

%Since we have shown that $Y^n$ is positive Harris recurrent,
By~\cite{DaiMeyn} the following ergodic property holds for every
measurable $f$ on $\mathrm{X}$ with $\pi(|f|) < \infty$,
\[
\lim_{t \rightarrow\infty} \frac{1}{t} \int_0^t
f\bigl(X^{n}(s)\bigr)\,ds = \pi(f) \qquad\mbox{$\Pb_x$-a.s.
for each $x \in\mathrm{X}$},
\]
where $\pi$ is the unique invariant distribution of $X^n$. Assigning
the function
$
f(x)\triangleq M^{-1}\dot{T}^{\mathbf{x}}(0) %{
$
to be the instantaneous service rates when the process is in state $x$
(recall that we assumed the service rates are a function of the state
in Section~\ref{secpreliminaries}), we have
%
%e45 #&#
%
\begin{equation}
\label{eqRcalconv} \lim_{t \rightarrow\infty} \frac{1}{t}\int
_0^t f\bigl(X^{\mathbf
{x}}(s)\bigr)\,ds =
\lim_{t \rightarrow\infty}\frac{1}{t} M^{-1}\widetilde{T}{}^{\mathbf
{x}}(t)
= \mathcal{R} \qquad\mbox{a.s.}
\end{equation}
for some constant vector $\mathcal{R}$.

%%%%%%%%%%%%%%%%%%%
%
% AUGMENT CHAIN
%
%%%%%%%%%%%%%%%%%%%%%%%%%

Consider the random variable $ \Ni\triangleq\liminf_{t\rightarrow
\infty} \frac{t}{N(t)}. $ The random variable $\Ni$ is a $\Pb_\pi$
invariant random variable, and therefore is a constant. Moreover by
(\ref{eqfatou}),
$
\mbox{ $\Ni\leq n t_0'$.}
$
A more detailed explanation of this argument is provided in \cite
{MusacchioPhD}.

We observe that for any sample path the following inequalities hold:
%
%e46 #&#
%
\begin{equation}
\label{eqsandwich} \frac{t}{N(t)} \frac{M^{-1} {T}^{\mathbf{x}}(t)}{t}
\leq\frac{ \sum_{j=0}^{N(t)} \rho_j }{N(t)}
\leq\frac{t}{N(t)} \frac{\sigma_{N(t)+1}}{t} \frac{M^{-1}{T}^{\mathbf
{x}}(\sigma_{N(t)+1})}{\sigma_{N(t)+1}}.
\end{equation}
Taking the $\liminf_{t\rightarrow\infty}$ of both sides and using
(\ref
{eqRcalconv}), we have that
%
%e47 #&#
%
\begin{equation}
\label{eqgreat} \liminf_{t \rightarrow\infty} \frac{ \sum_{j=0}^{N(t)}
\rho_j }{N(t)}= \Ni\mathcal{R}
\qquad\mbox{a.s.}
\end{equation}
We note that
$ {{T}^{\mathbf{x}}(\sigma_{N(t)+1})}/ {\sigma_{N(t)+1}} \leq
\mathbf{1}$
where $\mathbf{1}$ is a column vector of 1's of appropriate dimension.
This fact combined with (\ref{eqsandwich}) yields that for each $i>0$,
\[
\inf_{k\geq i} \frac{\sum_{j=1}^k {\rho}_j}{i} \leq\liminf_{t\rightarrow
\infty}
\frac{t}{N(t)} M^{-1}\mathbf{1} \leq n t_0'
M^{-1}\mathbf{1}. %t_0(1+\phi)M^{-1}\mathbf{1}.
\]
Thus the random variables
$ \{ \inf_{k\geq i} i^{-1} \sum_{j=1}^k {\rho}_j\dvtx  i>0
\}
$
are dominated by a constant. Consequently, %the dominated convergence
%theorem applied to (\ref{eqgreat}) yields,
$
\liminf_{i \rightarrow\infty} \E[ {\sum_{j=1}^i {\rho
}_j}/{i} ] =
\Ni\mathcal{R}
$
by the dominated convergence theorem. Also\vspace*{1pt} for each $i>0$,
$
\E[ {\sum_{j=1}^i {\rho}_j}/{i} ] \geq R n t_0(1-\gamma
)(1-\zeta)
$
by~(\ref{eqErho}).
Thus,
$ %\[
\Ni\mathcal{R} \geq Rn t_0(1-\gamma)(1-\zeta).
$ %\]
Substituting (\ref{eqfatou}) we have that
$ %\[
\mathcal{R} \geq(1-\gamma)(1-\zeta)\frac{t_0}{t_0'} R. %\frac{(1-
$ %\]
This implies % by (\ref{eqvarphi}) - (\ref{eqRcalconv}) implies
\[
\lim_{t \rightarrow\infty} \frac{1}{t} M^{-1} T^{\mathbf{x}}(t)
\geq\frac{(1-\gamma)(1-\zeta)}{1+ ({\zeta+b})/({1-\delta})} R
\qquad\mbox{a.s.}
\]
Recall $\gamma$, $\zeta$, $b$ and $\delta$ may be chosen arbitrarily
small, so long as $n$ is chosen large enough according to (\ref
{eqnchoice}). Thus, for any $\epsilon>0$ there exisits an $n$ such that
\[
\lim_{t \rightarrow\infty} \frac{1}{t} M^{-1}T^{\mathbf{x}}(t) \geq
(1-\epsilon)R \qquad\mbox{a.s.}
\]
%
%for any initial condition $x$.
By the strong law of large numbers for renewal processes~\cite{Durrett},
$
% \lim_{t\rightarrow\infty}
\frac{1}{t} S_k^{\mathbf{x}}(t) \rightarrow
m_k
$
a.s.
Thus by (\ref{eqdepartures}),
$\lim_{t\rightarrow\infty}
\frac{1}{t} D^{\mathbf{x}}(t) \geq(1-\epsilon)R$ a.s.
\end{pf}

%s4 #&#
\section{Analysis of switch example}
\label{secswitchexample}

In this section we apply the results of the preceding section to the
example introduced in Section~\ref{secbigexample}. Recall that this example
resembles a 2-input 2-output switch and has 3 flows and is illustrated
by Figure~\ref{figexample2}. As we discussed in Section \ref
{secbigexample}, the max-min fair share rate allocation would be that
all three flows achieve rates of $0.5$, so we set $R=[0.5, 0.5, 0.5]^T$
to be the vector of desired rates.

\begin{figure}

\includegraphics{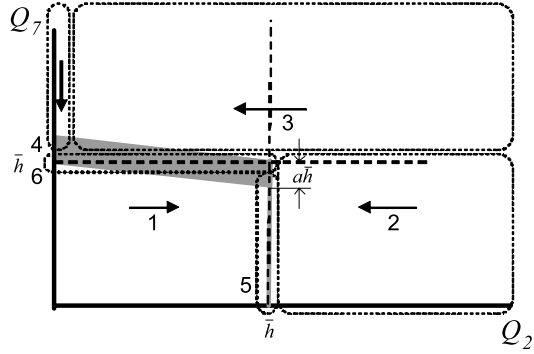}

\caption{The evolution of $(\bar{Q}_2(t),\bar{Q}_7(t))$. The shaded area
indicates the set $\bar{h} \eq$.}
\label{figphaseport}
\end{figure}

To fit the framework we have developed, we must show that the fluid
model with thresholds $\bar{h}$ is drawn to a set $\bar{h}\eq$, and
that the fluid model rates while in $\bar{h}\eq$ are $R$.
%To apply the results of the preceding section, we must show that the
%fluid
%model of the system, described by equations (\ref{eqUV})-(
%admits trajectories that go to a trapping set from any initial
%condition.
%Furthermore, the trajectories should reach the set in a
%time not more than a constant times the initial condition's distance
%from the set.
Intuition suggests that the dynamics of the fluid model should evolve
in the following way:
\begin{itemize}
\item One of the queues flow 2 passes through (either queue 2 or 7)
reaches threshold and ``chatters'' there. The other queue can be
anywhere at or below its threshold. By ``chatters'' we mean that it
alternately goes a tiny amount above and below. However, if the
differential inclusions of the fluid model are such that: (i) the queue
grows whenever below threshold or (ii) shrinks when above, then a fluid
model trajectory would go to threshold and stay there.
\item Queue 1 fills to threshold, ``chatters'' there, limiting flow 1's
ultimate rate.
\item Queue 7 fills to threshold, ``chatters'' there, limiting flow 3's
ultimate rate.
\item Other queues are not ``bottlenecks'' and should empty.
\end{itemize}
%
%Intuition
%suggests that one of the queues that flow 2 passes through should
%reach its
%threshold and chatter there. Similarly, flow 1's rate should
%ultimately be
%limited by queue 1, and thus queue 1 should fill to its threshold and
%chatter
%there. The other queue that flow 1 passes through should be nearly
%empty
%because that station has more capacity than the arrival rate of the
%one flow it
%handles. By the same reasoning, flow 3 should be limited by queue 7,
%and thus
%queue 7 should fill to its threshold and chatter there, while queue 3
%should be
%nearly empty.
This above intuition suggests that the fluid model is drawn to the set
$\bar{h} \tilde{\eq}$ where $\tilde{\eq}$ is given by
\begin{eqnarray*}
\tilde{\eq} &\triangleq& \bigl\{\bar{X}\dvtx\bar{Q}_1 =
\bar{Q}_8 = 1, \bar{Q}_3=\bar{Q}_4=
\bar{Q}_5=\bar{Q}_6= 0,
\\
&&\hspace*{4pt}(\bar{Q}_2, \bar{Q}_7) \in\bigl\{ [0,1] \times1
\bigr\} \cup\bigl\{ 1 \times[0,1] \bigr\}, \bar{U}=0, \bar{V}=0, \bar
{H}=0\bigr
\}.
\end{eqnarray*}
%
%In this definition for $\bar{h} \tilde{\eq}$, one of the queues that
%flow 2 passes
%through must be at its threshold, while the other queue must be at or
%below its
%threshold.
As it will turn out, the most critical part of the analysis of this
example's fluid model is to show that the queues flow 2 passes through, queues
2 and 7, go to values in $\bar{h} \tilde{\eq}$ in a time not more than
a constant times their initial values.
Intuition suggests that
after a ``settling down'' period flow 1's rate through queue 1, as well
as flow
3's rate through queue 8, settles to $0.5$. After flow 1 and flow 3's
rates settle, the dynamics of
%
%t1 #&#
%
\begin{table}[b]
\caption{Dynamics of $(\bar{Q}_2(t),\bar{Q}_7(t))$, after flows 1
and 3
settle to their ultimate rates of 0.5. The rows numbers correspond to
the regions labeled in the phase portrait diagram of Figure \protect
\ref
{figphaseport}}% Dependence of
%$(\bd{Q}_2(t),\bd{Q}_7(t))$ on $(\bar{Q}_2(t),\bar{Q}_7(t))$. The
%column
%labeled ``Time to Equilibrium'' indicated the amount of time to reach
%the set
%$\eq$. The column labeled ``Time to Distance Ratio'' contains the
%value of the
%maximum ratio in the indicated region of the amount of time to reach
%the set
%$\eq$ divided by the starting $L^1$ distance to $\eq$.}
\label{tbphaseport}
\begin{tabular*}{\tablewidth}{@{\extracolsep{\fill}}lcccclc@{}}
\hline
& $\bolds{\bar{Q}_2}$ & $\bolds{\bar{Q}_7}$
& $\bolds{\bd{Q}_2}$ & $\bolds{\bd{Q}_7}$ & \textbf{Time to}
$\bolds{\bar{h} \eq}$ &
$\bolds{\frac{\mathrm{Time}\ \mathrm{to}\ \bar{h} \eq}{ \llVert \bar{X} \rrVert_{
\bar{h} \eq} }}$\\
\hline
1 & $[0,\bar{h})$ & $[0,\bar{h})$ & \hphantom{$-$}0.1 & 0 &
if $|\bar{Q}_7-\bar{h}| < a\bar{h}$ then & $\frac{10}{a}$\\[2pt]
&&&&& $\frac{10}{a}|\bar{Q}_7-\bar{h}| $\\[2pt]
&&&&&
if $|\bar{Q}_7-\bar{h}| \geq a \bar{h}$ then\\[2pt]
&&&&& $10|\bar{Q}_2-\bar{h}|$
\\
%
% $\cases{
% 10\frac{\bar{h}}{a}|\bar{Q}_7-\bar{h}| & \mbox{if $|\bar{Q}_7-
% 10|\bar{Q}_2-\bar{h}| & \mbox{if $|\bar{Q}_7-\bar{h}| \geq a$}
% }$
[4pt]
2 & $(\bar{h},\infty)$ & $[0,\bar{h})$ & $-$0.5 & 0 & $2|\bar
{Q}_2-\bar
{h}|$ & 2 \\[4pt]
3 & $(0,\infty)$ & $(\bar{h},\infty)$ & $-$0.5 & 0 &
if $|\bar{Q}_7-\bar{h}| < a\bar{h}$ then & $\frac{2}{a}$\\[2pt]
&&&&& $\frac{2}{a}|\bar{Q}_7-\bar{h}|$ \\
[4pt]
&&&&& if $|\bar{Q}_7-\bar{h}| \geq a\bar{h}$ then\\
&&&&& $2|\bar{Q}_2-\bar{h}| + 2|\bar{Q}_7-\bar{h} - a\bar{h}|$
\\
%
% $\cases{
% 2\frac{h}{a}|\bar{Q}_7-\bar{h}| & \mbox{if $|\bar{Q}_7-\bar{h}|< a$}
% 2\left[
% \left[\begin{matrix} \bar{Q}_2 \\ \bar{Q}_7 \end{matrix} \right] -
% \left[\begin{matrix} \bar{h} \\ \bar{h} \end{matrix}\right]
% \right]
% 2\left[
% \begin{array}{l} |\bar{Q}_2-\bar{h}| \\ + |\bar{Q}_7-\bar{h}|
% \right]
% & \mbox{if $|\bar{Q}_7-\bar{h}| \geq a$}
% }$
[4pt]
4 & 0 & $(\bar{h},\infty)$ & 0 & $-$0.5 & $2|\bar{Q}_7 - \bar{h}-a|$ &
2\\
[4pt]
5 & $\bar{h}$ & $[0,\bar{h}]$ & $[-0.5,0.1]$ &
$[-0.5,0]$ & \multicolumn{1}{c}{0} & N$/$A\\
[4pt]
6 & $[0,\bar{h}]$ & $\bar{h}$ & $[-0.5,0.1]$ &
$[-0.5,0]$ & \multicolumn{1}{c}{0} & N$/$A\\
\hline
\end{tabular*}\vspace*{-4pt}
\end{table}
$(\bar{Q}_2(t), \bar{Q}_7(t))$, the queues of flow 2, follow the
relations outlined by
Table~\ref{tbphaseport} and illustrated by Figure \ref
{figphaseport}. The
entries of Table~\ref{tbphaseport} are easily derived by using the
observations that:
\begin{itemize}
\item The arrival rate to queue 2 is 0.6 when queue 2 and queue 7 are below
threshold while the arrival rate to queue 2 is 0 when one of these
queues is
above threshold.
\item The departure rate from either queue 2 or queue 7 is 0.5
whenever the queue is nonempty or has sufficient arrivals to maintain this
departure rate. (This relies on our assumption that the flow rates through
queues 1 and 8 have ``settled down'' to $0.5$.)
\end{itemize}
Figure~\ref{figphaseport} is a\vspace*{2pt} vector flow diagram, showing the
dependence of $(\bd{Q}_2(\cdot), \bd{Q}_7(\cdot))$ on $(\bar
{Q}_2(\cdot
), \bar{Q}_7(\cdot))$.
It is evident from the diagram that the time to reach the set
\[
\bigl\{ [0,\bar{h}] \times\bar{h} \bigr\} \cup\bigl\{ \bar{h} \times
[0,\bar{h}]
\bigr\},
\]
which is the projection of $\bar{h} \tilde{\eq}$ on to the subspace on
which $(\bar{Q}_2(\cdot), \bar{Q}_7(\cdot))$ takes values, is not
always less than or equal to a constant times the initial condition's
distance from this set. Consider an
initial condition of $(\frac{\bar{h}}{2}, \bar{h} + \epsilon)$.
This initial
condition is only a distance of $\epsilon$ from $\tilde{\eq}$, but the
time it
takes to reach the set $\tilde{\eq}$ is $\bar{h} + \frac
{1}{2}\epsilon
$. (Note
%
%f5 #&#
%
that we will use the $L^1$ norm throughout this section.) This is the
same phenomenon we observed in
the example in the \hyperref[sintro]{Introduction} of the paper. There, as here, we can
fix the
problem by slightly enlarging the set $\tilde{\eq}$ to a new set $\eq$
so that
the set is reached in a time not more than a constant times the initial
condition's starting
distance from the set. %However, we must not enlarge must not enlarge
%the set too much, because we want the worst case flow rates for states
%in the set to be close to the desired rates.
%However, the worst case flow rates for states in $\eq$
%will ultimately be the rates that we show the original stochastic
%system
%reaches. If we make $\eq$ too large so as to include states whose flow
%rates
%are low, the lower bound on the asymptotic flow rates of the
%stochastic system
%will not be as strong.
To this end, we define $\eq$
according to \label{eqdef}
\begin{eqnarray*}
\eq&\triangleq& \bigl\{\bar{X}\dvtx\bar{Q}_1 =
\bar{Q}_8 = 1, \bar{Q}_3=\bar{Q}_4=
\bar{Q}_5=\bar{Q}_6= 0,
\\
&&\hspace*{4.7pt}( \bar{Q}_2, \bar{Q}_7) \in\bigl\{(\chi,\psi)\dvtx
\chi\in[0,1], \psi\in\bigl[1 - a\chi, 1 + a( 1-\chi)\bigr]\bigr\}
\\
&&\hspace*{112.3pt}{}\cup\bigl\{ 1 \times[0,1] \bigr\}, \bar{U}=0, \bar{V}=0, \bar
{H}=0\bigr\}.
\end{eqnarray*}
Here $a$ is an arbitrary positive constant that should be less than
$1$. The projection of this set onto the subspace spanned by
$(Q_2,Q_7)$ is shown as the shaded area in Figure~\ref{figphaseport}.
With this definition, one can show that the set $\bar{h} \eq$ is
reached in a time not more than a constant times the initial distance
from $\bar{h} \eq$. The time to
reach $\bar{h} \eq$, along with the maximum ratio of the time to reach
$\bar{h} \eq$ divided by
initial distance to $\bar{h}\eq$ are shown in Table~\ref{tbphaseport}.
% Note that we
%are slightly abusing the term ``equilibrium'' slightly by calling $
%of equilibria because the fluid model's state can continue to change
%when
%starting from a state in $\eq$.

We are now ready to formalize the intuition we have outlined in the preceding
paragraphs. We begin by stating a lemma that the system settles down so that
the behavior flow 2's queues are as described by Table \ref
{tbphaseport} after
a time $\tau_{sd}$ (mnemonic for ``settle down'') that is in proportion
to the
initial condition.

%le3 #&#
%
\begin{lemma}
\label{lmsettledown} There exists a time $\tau_{sd}$ proportional to the
initial condition as described by the relation
\[
\tau_{sd} = t_{01} \bigl\llVert\bar{X}(0) \bigr
\rrVert_{ \bar{h} \eq}
\]
for some positive $t_{01}$ such that for all
regular points $t\geq\tau_{sd}$:
\begin{itemize}
\item The value of $(\bd{Q}_2(t), \bd{Q}_7(t))$ is determined by the
value of
$(\bar{Q}_2(t), \bar{Q}_7(t))$ as specified by Table
\ref{tbphaseport}.\vspace*{1pt}

\item$\bar{Q}_3(t)=\bar{Q}_5(t)=\bar{Q}_6(t)=0$, and $\bar{Q}_8(t) =
\bar{h}$.

\item The time to reach the set $\bar{h} \eq$, as well as the maximum ratio
between this time and the distance of $(\bar{Q}_2(\tau_{sd}),
\bar{Q}_7(\tau_{sd}))$ from $\bar{h} \eq|_{\bar{Q}_2,\bar{Q}_7}$ in
any of the regions
1 through 4 is as specified in Table~\ref{tbphaseport}. [Here $\bar
{h} \eq|_{\bar{Q}_2,\bar{Q}_7}$ denotes the projection of the set
$\bar
{h} \eq$ onto the space on which $(\bar{Q}_2,\bar{Q}_7)$ takes values.]
\end{itemize}
\end{lemma}

Lemma~\ref{lmsettledown} is proved by using relations (\ref
{eqUV})--(\ref{eqRRdrain1}) that describe the evolution of a fluid
model trajectory. The proof is straightforward but slightly lengthy
because it requires analysis for each entry in Table \ref
{tbphaseport}. We therefore omit this proof.

% result is simple though slightly tedious to prove. We therefore omit
%the proof.
%intuitive but tedious to prove. We
%therefore give the proof in an appendix.
We now state and prove the principal
result of this section.
%
%th5 #&#
%
\begin{theorem}
For any $\epsilon>0$, there exists an $n_c>0$ such that if the
discarding thresholds of the stochastic system
in Example 2 are set to $n h$, $n \geq n_c$, then
\[
\lim_{t \rightarrow\infty} \frac{D(t)}{t} \geq(1-\epsilon)
\frac
{1}{2}\mathbf{1} \qquad\mbox{a.s.},
\]
where $\mathbf{1}$ is a vector of ones of dimension $K$.
\end{theorem}
\begin{pf}
By Lemma~\ref{lmsettledown} the dynamics of the state variables
$(\bar{Q}_2(t), \bar{Q}_7(t))$ of the fluid model trajectory evolve according
to Table~\ref{tbphaseport} after a time $\tau_{sd} = t_{01} \llVert
\bar{X}(0) \rrVert_{ \bar{h} \eq}$. From Table \ref
{tbphaseport}, the $\bar{Q}_2$ and
$\bar{Q}_7$
components of the fluid model trajectory reach values in the set $\bar
{h}\eq$'s
projection in, at most, an additional $\frac{10}{a}
\llVert \bar{X}(\tau_{sd}) \rrVert_{ \bar{h} \eq}$ time units.
%(We say in the ``trapping
%set's projection'' to emphasize that some components of the state
%vector have
%reached values consistent with $\eq$, but that not all components have
%yet done
%so.)
Because the total arrival rate into the system is less than or equal to
$1.8$,
\[
\bigl\llVert\bar{X}(\tau_{sd}) \bigr\rrVert_{ \bar{h} \eq}
\leq(2t_{01} + 1)\bigl\llVert\bar{X}(0) \bigr\rrVert_{ \bar{h} \eq}.
\]

Thus after a time $t_{02}\llVert \bar{X}(0) \rrVert_{ \bar{h} \eq
}$, where
$ t_{02} = \frac{10}{a} (2t_{01} + 1) + 1, $
all queues but queue 1 have been shown to reach values in the
projection of the
set $\bar{h}\eq$. By Lem\-ma~\ref{lmsettledown}, queue 5 is empty, so either:
queue 1 is above threshold, in which case discarding is on and it will reach
threshold in $ 2 (\bar{Q}_1 ( t_{02}\llVert \bar{X}(0)
\rrVert_{ \bar{h} \eq}
) - \bar{h} ) $ time units, or queue 1 is below threshold
in which
case it will reach threshold in $ 10 (\bar{Q}_1 (
t_{02}\llVert \bar{X}(0) \rrVert_{ \bar{h} \eq} ) - \bar
{h} ) $ time units. Once $\bar
{Q}_1(t)$ reaches
threshold $\bar{h}$, it remains there by the following reasoning. If
queue 1
were to move some positive amount $\epsilon$ above $\bar{h}$, the discarding
would have turned on before the queue grew to $\epsilon$ and prevented
it from
getting there. Similarly, if queue 1 were to move some positive amount
$\epsilon$ below $\bar{h}$, the discarding would have turned off before
the queue
receded by $\epsilon$, and prevented the queue from receding that
much. Very
loosely, we can bound the rate of growth of queue 1 before time
$t_{02}\llVert \bar{X}(0) \rrVert_{ \bar{h} \eq} $ by
\[
\bar{Q}_1 \bigl( t_{02}\bigl\llVert\bar{X}(0) \bigr
\rrVert_{ \bar{h} \eq} \bigr) \leq1.6 \bigl\llVert\bar{X}(0) \bigr
\rrVert_{ \bar{h}}.
\]

Thus after a time of length $t_0 \llVert \bar{X}(0) \rrVert_{ \bar
{h} \eq}$, where
$t_0$ is
given by
$ t_0 = \frac{10}{a} (2t_{01} + 1) + 17$,
all fluid model trajectories will have reached the set $\bar{h}\eq$.
The departure
rates for all three flows, as well as the departure rates for each
class associated with each flow, are easily seen to be $0.5$ when the
fluid model's
state is in $\bar{h}\eq$ and threshold $\bar{h}>0$. Thus, by
Theorem~\ref{thbigtheorem} we have that the asymptotic flow rates approach
$0.5$.
\end{pf}

%s5 #&#
\section{Conclusion}
\label{sconclusion}

In this work we have shown how the analysis of the flow rates of a
stochastic network with a particular flow control scheme may be
reduced to an analysis of a fluid model. While we have focused on a
particular flow control scheme, the same analysis could
be carried out for many other control schemes. The key feature that
enabled our approach was that our control scheme has a free
parameter, $n$, which when increased makes the system look more and
more like a deterministic fluid system. We have demonstrated how to
use the theory developed in this paper to analyze an example network
resembling a 2-input, 2-output switch.

\begin{appendix}\label{app}
\section*{Appendix}

Before proving Theorem~\ref{thsubseqconv}, we state and prove a
number of lemmas.
Lemma~\ref{lemFSLLN} is a functional form of the strong law of large
numbers for renewal processes, and is\vadjust{\goodbreak} taken from~\cite{Dai95}.
Lemma~\ref{lemthin} is a new result showing that the thinned arrivals
(the customers that make it beyond the discarding point) converge to a
fluid limit along a subsequence.
Lemma~\ref{lemuniformint} is a result taken from~\cite{Dai95} showing
that the residual initial arrival and service times decline to zero at
rate $1$ in the fluid limit. The lemma also shows that the sequence of
functions we use to take the fluid limit are uniformly integrable.

Also the lemmas will make use of fluid limits that have well-defined
limiting residual interarrival and service times, as defined by the
following property.
%
%pr2 #&#
%
\begin{property}
\label{propb}
$\{ ({{\mathbf{x}}_{j}},a_j)\}$ is a sequence for which $
\frac{U^{{\mathbf{x}}_{j}}(0)}{a_j} \rightarrow\bar
{U}(0)$, $\frac{V^{{\mathbf{x}}_{j}}(0)}{a_j} \rightarrow
\bar{V}(0)$, for some $\bar{U}(0) \geq0$
and $\bar{V}(0) \geq0$.
\end{property}

%%%%%%%%%%%%%%%%%%%%%%%%%%%%%%%%%%%%%%
%
% Lemma 1
%
%%%%%%%%%%%%%%%%%%%%%%%%%%%%%%%%%%%%%
%
%le4 #&#
%
\begin{lemma}[(Dai, Lemma 4.2 of~\cite{Dai95})]
\label{lemFSLLN} Suppose that $\{ ({{\mathbf{x}}_{j}},a_j
)\}$ is a
sequence satisfying \label{referback2} Properties~\ref{propa} and
\ref{propb}
\ifthenelse{\equal{\pageref{propa}}{\pageref{propb}}}{(on page
\pageref{propa}) } {(see pages \pageref{propa} and \pageref{propb})}.
Then for almost all $\omega$,
\begin{eqnarray*}
\frac{E_f^{\mathbf{x}_{j}}(a_{j}t)}{a_j} &\rightarrow&\alpha_f\bigl(t-\bar
{U}_f(0)\bigr)^+ \qquad\mbox{u.o.c.},
\\
\frac{S_k^{\mathbf
{x}_{j}}(a_{j}t)}{a_j} &\rightarrow&
\mu_k\bigl(t-\bar{V}_k(0)\bigr)^+ \qquad\mbox{u.o.c.}
\end{eqnarray*}
\end{lemma}
\begin{pf}
See Lemma 4.2 of Dai~\cite{Dai95}. The result is an instance of the
strong law of large numbers for renewal processes~\cite{Durrett}.
\end{pf}
%
%%%%%%%%%%%%%%%%%%%%%%%%%%%%%%%%%%%%%%
%
% Lemma 2
%
%%%%%%%%%%%%%%%%%%%%%%%%%%%%%%%%%%%%%
%
%le5 #&#
%
\begin{lemma}[(Thinned arrival convergence)]
\label{lemthin} Suppose that $\{ ({{\mathbf{x}}_{j}},a_{j
})\}$ is a
sequence satisfying Properties~\ref{propa} and~\ref{propb}. Then for
almost all $\omega$, there exists a subsequence $\{ ({{\mathbf
{x}}_{m}},a_{
m})\}\subseteq\{ ({{\mathbf{x}}_{j}},a_{j})\}$ such that
\[
{\Lambda^{\mathbf{x}_{m}}(a_{m}t)}/{a_{m}} \rightarrow
\bar{\Lambda}(t) \qquad\mbox{u.o.c.},
\]
where $\bar{\Lambda}(t)$ is some Lipschitz continuous process
with, for all regular $t\geq0$,
%
%e48 #&#
%
\begin{equation}
\label{eqGlimit2lem} \bd{\Lambda}_f(t) \leq\alpha_f
\qquad\mbox{for each flow $f$.}
\end{equation}
\end{lemma}
\begin{pf}
By Lemma~\ref{lemFSLLN},
%
%e49 #&#
%
\begin{equation}
\label{eqexoguoc} {E_f^{\mathbf{x}_{j}}(a_{j}t)}/{a_{j}}
\rightarrow\alpha_f\bigl( t -\bar{U}_f(0)\bigr)^+
\qquad\mbox{u.o.c.}
\end{equation}
for each flow $f$. For notational convenience in the development that
follows, we define
%
%e50 #&#
%
\begin{equation}
\label{eqDelta} \bar{E}_f(t) \triangleq\alpha_f\bigl(
t - \bar{U}_f(0)\bigr)^+,\qquad \Delta_j(t)
\triangleq{E^{\mathbf{x}_{j}}(a_{j}t)}/{a_{j}} - \bar{E}(t).
\end{equation}
Pick a compact time interval $[s_0,s_1]$. Since the number of admitted
customers is not greater than the number that arrive,
%
%e51 #&#
%
\begin{equation}
\label{eqbigmess} \frac{1}{a_{j}} \bigl[ \Lambda^{\mathbf{x}_{j}}
\bigl(a_{j}(t+\varepsilon)\bigr) - \Lambda^{\mathbf{x}_{j}}(a_{j}t)
\bigr]\leq\frac{1}{a_{j}} \bigl[E^{\mathbf{x}_{j}}\bigl(a_{j}(t+
\varepsilon)\bigr) - E^{\mathbf{x}_{j }}(a_{j }t) \bigr]
\end{equation}
for any positive $\varepsilon\leq s_1-s_0$ and $t\dvtx  s_0\leq t \leq
s_1-\varepsilon$. Adding $-\Delta_j(t+ \varepsilon)$ and
$\Delta_j(t)$ to both sides and substituting (\ref{eqDelta}), we have
\[
\frac{\Lambda^{\mathbf{x}_{j}}(a_{j}(t+\varepsilon))}{a_{j}} - \Delta_j
(t+ \varepsilon) - \biggl[
\frac{\Lambda^{\mathbf{x}_{j}}(a_{j}t)}{a_{j}} - \Delta_j (t ) \biggr]
\leq\bar{E}(t+\varepsilon)
- \bar{E}(t) \leq\varepsilon\alpha.
\]

%*****Question: I do not get this argument. On the left of the
%inequality, you have random quantities. The upper bound is
%deterministic. This makes sense in the limit, but why should this
%always hold?*****

Define the family of functions
\[
\frak{L}_{j}(s_0,t):= \sup_{s\in[s_0,t]} \biggl[
\frac{\Lambda^{\mathbf{x}_{j}}(a_{j}s)}{a_{j}} - \Delta_j(s ) \biggr]
\]
for $t\in[s_0,s_1]$.
Because the argument of the $\sup$ function is a vector, $\sup$ is
taken component-wise. Note that for any $(t,\varepsilon)$ with $t \in
[s_0,s_1-\varepsilon]$,
\[
\frak{L}_{j}(s_0,t+\varepsilon) = \frak{L}_{j}(s_0,t)
\vee\frak{L}_{j}(t,t+\varepsilon)
\]
and
\[
\frak{L}_{j}(t,t+\varepsilon) \leq\varepsilon
\alpha+ \frak{L}_{j}(t,t) \leq\varepsilon\alpha+ \frak{L}_{j}(s_0,t).
\]
Thus $ \frak{L}_{j}(s_0,t+\varepsilon) - \frak{L}_{j}(s_0,t)
\leq\varepsilon\alpha$ and clearly $\frak{L}_{j
}(s_0,t+\varepsilon
) - \frak{L}_{j}(s_0,t) \geq0$ because $\frak{L}_{j
}(s_0,\cdot
)$ is monotone. Hence the functions $\frak{L}_{j}(s_0,\cdot)$ are
equicontinuous and individually Lipschitz continuous. Thus, by Arzela's
theorem, there exists a further subsequence $\{({{\mathbf
{x}}_{m}},a_{m})\}\subseteq
\{({{\mathbf{x}}_{j}},a_{j})\}$ such that
\[
\frak{L}_{m}(s_0,t) \rightarrow\bar{\Lambda}(t)
% \sup_{s\in[s_0,t]}\left[ \icsft{\Lambda}{l}{s} - \Delta^j(s )
% \rightarrow\bar{\Lambda}(t)
\]
uniformly on the compact set $t\in[s_0,s_1]$ for some
monotone-nondecreasing, Lipschitz-continuous process $\bar{\Lambda
}(t)$. But by (\ref{eqexoguoc}), $\Delta_j(t ) \rightarrow0$
uniformly on compact sets. Because of this and the fact that $\Lambda
^{\mathbf{x}_j}(a_{j} s)/a_{j}$ is monotone in
$s$, it follows
that %the maximizing values of each $\sup$ term in the definition of
$\frak{L}_{j}(s_0,t)$ approaches $\Lambda^{\mathbf{x}_j}(a_{j}
t)/a_{j}$ as $j\rightarrow\infty$. Thus
\[
\sup_{s\in[s_0,t]} \biggl[ \frac{\Lambda^{\mathbf
{x}_{j}}(a_{j}s)}{a_{j}} - \Delta_j (s )
\biggr] \rightarrow\frac{\Lambda^{\mathbf{x}_{j}}(a_{j}t)}{a_{j}}
\rightarrow\bar{\Lambda}(t).
\]
Because the choice of $[s_0,s_1]$ was arbitrary, we have $\Lambda
^{\mathbf{x}_m}(a_{m} s)/a_{m} \rightarrow\bar
{\Lambda}(t)$ u.o.c.
Furthermore, (\ref{eqexoguoc}) and (\ref{eqbigmess}) imply that
$\bar
{\Lambda}(t)$ satisfies (\ref{eqGlimit2lem}).
\end{pf}

%%%%%%%%%%%%%%%%%%%%%%%%%%%%%%%%%%%%%%
%
% Lemma 3
%
%%%%%%%%%%%%%%%%%%%%%%%%%%%%%%%%%%%%%
%
%le6 #&#
%
\begin{lemma}[(Lemmas 4.3 and 4.5 of Dai~\cite{Dai95})]
\label{lemuniformint} %\textbf{u.i. should be for a sequence}
Suppose that $\{({{\mathbf{x}}_{j}},a_{j})\}$ is a sequence satisfying
Properties \ref
{propa} and~\ref{propb}. Then almost surely
\begin{eqnarray*}
\lim_{j\rightarrow\infty} \frac{U_f^{\mathbf{x}_{j}}(a_{j}t)}{a_j}
&=&\bigl(\bar{U}_f(0) - t
\bigr)^+ \qquad\mbox{u.o.c.},
\\
\lim_{j\rightarrow\infty} \frac{V_k^{\mathbf
{x}_{j}}(a_{j}t)}{a_j} &=&\bigl(\bar
{V}_k(0) - t\bigr)^+ \qquad\mbox{u.o.c.}
\end{eqnarray*}
Also, for each fixed $t\geq0$, the sets of functions
\begin{eqnarray*}
&\displaystyle \bigl\{ U^{{\mathbf{x}}_{j}}(a_jt) / a_j\dvtx  a_j
\geq1 \bigr\},\qquad \bigl\{ V^{{\mathbf{x}}_{j}}(a_jt) /
a_j\dvtx  a_j\geq1 \bigr\},&
\\
&\displaystyle \bigl\{ Q^{{\mathbf{x}}_{j}}(a_jt) / a_j\dvtx  a_j
\geq1 \bigr\}&
\end{eqnarray*}
are uniformly integrable.
\end{lemma}
\begin{pf} See Lemmas 4.3 and 4.5 of Dai~\cite{Dai95}.
\end{pf}

%%%%%%%%%%%%%%%%%%%%%%%%%%%%%%%%%%%%%%
%
% Lemma 4
%
%%%%%%%%%%%%%%%%%%%%%%%%%%%%%%%%%%%%%
We use the following lemma later to show that because all of the
systems we consider are work-conserving, the fluid limit must also be
work-conserving. In the lemma below, the notation $D_\R[0,\infty)$
denotes the space of right-continuous functions
on $\R_+$ having left limits on $(0,\infty)$, and endowed with the
Skorohod topology~\cite{EthierKurtz}. $C_\R[0,\infty)\subset D_\R
[0,\infty)$ is the subset of continuous paths.
%
%le7 #&#
%
\begin{lemma}[(Lemma 2.4 of Dai and Williams~\cite{DaiWilliams})]
\label{lemintegralconv}
Let $\{(z_j,\chi_j)\}$ be a sequence in $ D_\R[0,\infty)
\times C_\R[0,\infty)$. Assume that $\chi_j$ is nondecreasing and
$(z_j,\chi_j)$ converges to $(z,\chi) \in C_\R[0,\infty)
\times C_\R[0,\infty)$ u.o.c. Then for any bounded continuous
function~$f$,
\[
\int_0^t f\bigl(z_j(s)\bigr)\,d
\chi_j(s) \rightarrow\int_0^t f
\bigl(z(s)\bigr)\,d\chi(s) \qquad\mbox{u.o.c.}
\]
\end{lemma}
\begin{pf}
See Lemma 2.4 of Dai and Williams~\cite{DaiWilliams}.
\end{pf}

We are now ready to prove Theorem~\ref{thsubseqconv}.

%%%%%%%%%%%%%%%%%%%%%%%%%%
%
% PROOF OF THEOREM 1
%
%%%%%%%%%%%%%%%%%%%%%%%%
%
\begin{pf*}{Proof of Theorem~\ref{thsubseqconv}}
Before scaling space, the discarding thresholds for each $j$ are
$n_jh$. After scaling space by $a_j$, the scaled thresholds
are $n_jh/ a_j$. Property~\ref{propa} insures that $n_j
/a_j$ is upper bounded by a constant. Thus by the
Bolzano--Weierstrass theorem, there exists a subsequence
$\{({{\mathbf{x}}_{r}},a_{r})\} \subseteq\{
({{\mathbf{x}}_{j}},a_{j})\}$ for which
$n_rh/a_r\rightarrow\bar{h}$ for some $\bar{h} \geq0$.

Property~\ref{propa} insures that $\llVert x_r/ a_r \rrVert_{ {n_r
h}/{a_r} \eq} $ is upper bounded by a constant. Thus
$ \limsup\llVert x_r/ a_r \rrVert_{ \bar{h} \eq}$ is finite.
Consequently, there must be some further subsequence $\{({{\mathbf
{x}}_{u}},a_{u})\}
\subseteq\{({{\mathbf{x}}_{r}},a_{r})\}$ for which
$ {x_u}/{a_u} \rightarrow\bar{X}(0)$ for some finite $\bar{X}(0)$.

The hysteresis variables satisfy $H^{\mathbf{x}_{u}}(a_{u}t) /
a_u\rightarrow0$
u.o.c. because $H^{\mathbf{x}_{u}}(a_{u}t)$ is bounded by a constant
by its
definition. This fact along with the convergence of $X^{{\mathbf
{x}}_{u}}(0)/
a_u\rightarrow\bar{X}(0)$ allows us to use Lemma \ref
{lemuniformint} to conclude
$ U^{\mathbf{x}_{u}}(a_{u}t) /\break a_u\rightarrow\bar{U}(t)$ and
$V^{\mathbf{x}_{u}}(a_{u}t) /
a_u\rightarrow\bar{V}(t)$ u.o.c.
% \s{r}^{-1} \bbmatrix U\ics{r} \\ V\ics{r} \ebmatrix
% \rightarrow\bbmatrix\bar{U}(t) \\ \bar{V}(t) \ebmatrix
%The u.o.c. convergence of $\s{r}^{-1}U\ics{r} \rightarrow
where $\bar{U}(t)$ and $\bar{V}(t)$ satisfy (\ref{eqUV}).
%Thus gives us the second and third components of
%Expression (\ref{eqtrajconv2}).

%In the TFL case, the second part of (\ref{eqcase1Xhe}) follows
%from (\ref{eqcase1ic}). For the JFL case, (\ref{eqcase2ica})
%combined with (\ref{eqlitlim}) imply
%%that $(\const h |y_r|^{-1}
%%-\bar{h}) \rightarrow0$ and thus
%%$|\bar{X}(0)-\bar{h}\ul{e}|=\const$, giving us
%the second part of (\ref{eqcase2Xhe}).

The cumulative service time process $T^{\mathbf{x}_u}$ satisfies
%
%e52 #&#
%
\begin{equation}
\label{eqTHconv} \bigl[T^{\mathbf{x}_u}(\s{u} t) -T^{\mathbf{x}_u}(\s
{u} s)
\bigr] / a_u \leq(t-s).
\end{equation}
Thus by Arzela's theorem~\cite{munkres}, there exists a further
subsequence $\{({{\mathbf{x}}_{v}},a_{v})\} \subseteq\{({{\mathbf
{x}}_{u}},a_{u})\}$ for which
${T^{\mathbf{x}_v}(a_vt)}/{{a_v}} \rightarrow\bar{T}(t)$.
Property (\ref{eqTinc}) follows from (\ref{eqTnd}). Property~(\ref{eqInd})
implies $
{I^{\mathbf{x}_v}(a_vt)}/{a_v} \rightarrow\bar{I}(t)
$ u.o.c. where $\bar{I}(t)$ satisfies (\ref{eqI}).

By Lemma~\ref{lemFSLLN}, ${S_k^{\mathbf{x}_{v}}(a_{v}t)}/{a_{v}}
\rightarrow(\mu_k t -
\bar{V}_k(0))^+$ u.o.c. for each class $k$. This fact combined
with (\ref{eqdepartures}) and (\ref{eqTHconv}) gives
(\ref{eqDeparturelimit1}).

We have already shown that $X^{\mathbf{x}_v}(0) \rightarrow
\bar{X}(0)$,
therefore the $U^{\mathbf{x}_v}(0)$ and $V^{\mathbf{x}_v}(0)$
components of
$X^{\mathbf{x}_v}(0)$ converge to a limiting value. This
fact allows us to
invoke Lemma~\ref{lemthin} to\vadjust{\goodbreak} conclude that there is a subsequence
$\{({{\mathbf{x}}_{m}},a_{m})\} \subseteq\{({{\mathbf
{x}}_{v}},a_{v})\}$ for which $
{\Lambda^{\mathbf{x}_m}}(a_mt) / {a_m}
\rightarrow\bar
{\Lambda
}(t)$ u.o.c. for some Lipschitz continuous process $\bar{\Lambda}(t)$
satisfying (\ref{eqGlimit2}).

Lemma~\ref{lemFSLLN} combined with (\ref{eqarrivaldef}) gives us
${A_k^{\mathbf{x}_{m}}(a_{m}t)}/{a_{m}} \rightarrow\bar{A}_k(t)$
u.o.c. for each class $k$
where $\bar{A}_k(t)$ is defined
by (\ref{eqArrivallimit}). Furthermore, $\bar{A}_k(t)$ is Lipschitz
continuous
because it is equal to a linear combination of functions we have
already shown
to be Lipschitz continuous. Thus using (\ref{eqqevolution}) we have that
%
%e53 #&#
%
\begin{equation}
\label{eqQconverges} {Q^{\mathbf{x}_{m}}(a_{m}t)}/{a_{m}}
\rightarrow\bar{Q}(t) \qquad\mbox{u.o.c.},
\end{equation}
where $\bar{Q}(t)$ is
a Lipschitz continuous function given by (\ref{eqQbar}). %Thus we
%have the first component of Expression (\ref{eqtrajconv2}).
Property (\ref{eqQbig0limit}) follows easily from
(\ref{eqQbig0}).

The next few arguments are similar to the proof of Proposition 4.2
in~\cite{DaiMeyn}.
Suppose that $\bar{Q}_k(t)>\bar{h}$ for some $k \in\mathcal{C}(f)$. By
Lipschitz continuity of $\bar{Q}_k(t)$, there exists some small $\tau
>0$ such
that
$
\min_{t\leq s\leq t+\tau} \bar{Q}_k(s) > \bar{h}.
$
By the uniformity of the queue convergence in (\ref{eqQconverges}) and
that $n_mh/a_m\rightarrow\bar{h}$, there
%Because $a_j^{-1} Q_k\icst{j}{s} \rightarrow\bar{Q}_k(s)$
%u.o.c. and $ a_j^{-1} n_j h \rightarrow\bar{h}$ u.o.c., there
exists $m^*$ such that for all $m>m^*$, $Q_k^{\mathbf{x}_{m}}(a_{m}s)
> n_mh$ for all $s \in
[t,t+\tau]$. Thus, by (\ref{eqLambda}) one finds that
$
\Lambda_f^{\mathbf{x}_{m}}(a_{m}s)-\Lambda_f^{\mathbf
{x}_{m}}(a_{m}t)=0, \forall
s\in
[t,t+\tau].
$
Therefore,\vspace*{-2pt} it follows that
$
\bar{\Lambda}_f(s)-\bar{\Lambda}_f(t) = 0, \forall s\in
[t,t+\tau]
$
and consequently, $\bd{\Lambda}_f(t) = 0$, which is (\ref{eqHlimit}).

Suppose that $\bar{Q}_k(t)<\bar{h}$ for all $k \in
\mathcal{C}(f)$. First note that in this case $\bar{h}>0$. %and
%therefore (\ref{eqthreshconverges}) implies that $n_j \rightarrow
By the Lipschitz continuity of $\bar{Q}_k(t)$ for each
$k$, there exists some small $\tau>0$ such that
$
\max_{k\in\mathcal{C}(f)} \max_{s\in[t,t+\tau]} \bar{Q}_k(s) <
\bar{h}.
$
Because $n_m\rightarrow\infty$, the uniformity of the
convergence in
(\ref{eqQconverges}), and that $n_mh/a_m\rightarrow\bar
{h}$, there
exists $m'$ such that for all $m>m'$, $Q_k^{\mathbf{x}_{m}}(a_{m}s)<
n_mh $. Furthermore
there exists a $m^*\geq m'$ such that for all $m>m^*$
and $k \in
\mathcal{C}(f)$, $Q_k^{\mathbf{x}_{m}}(a_{m}s)< n_mh - o(n_m) h
\varsigma
$. Thus, by~(\ref{eqLambda}),
%
%%
%E_f\icst{\seqm}{t} \forall s\in[t,t+\tau]
%%
%}{
%
\[
\Lambda_f^{\mathbf{x}_{m}}(a_{m}s)-\Lambda_f^{\mathbf
{x}_{m}}(a_{m}t)=
E_f^{\mathbf{x}_{m}}(a_{m}s)- E_f^{\mathbf{x}_{m}}(a_{m}t)\qquad
\forall s\in[t,t+\tau]
\]
and consequently we have (\ref{eqGlimit}).

Suppose that for some class $k$, $\bar{Q}_k(t)>0$. By the
Lipschitz continuity of $\bar{Q}_k(t)$ there exists some small
$\tau>0$ such that
$
\min_{t\leq s \leq t+\tau} \bar{Q}_k(s) > 0.
$
Because of the uniformity of convergence in (\ref{eqQconverges})
there exists
$m^*$ such that for all $m>m^*$,
$
Q_k^{\mathbf{x}_{m}}(a_{m}s)> 0,\forall s \in[t,t+\tau].
$
By (\ref{eqRRdrainf}), for almost all $\omega$, and all classes
$l$ we have
%%
%w_k^{-1}[D_k(a_\seqm s)-D_k(a_\seqm t)] \geq\\
%w_l^{-1}[D_l(a_\seqm s)-D_l(a_\seqm t)] \forall s \in[t,t+\tau]
%%
%}
%
\[
w_k^{-1}\bigl[D_k(a_ms)-D_k(a_mt)
\bigr] \geq w_l^{-1}\bigl[D_l(a_ms)-D_l(a_mt)
\bigr]\qquad\forall s \in[t,t+\tau],
\]
and thus we have (\ref{eqRRdrain2}).

If $\bar{Q}_l(t)>0$ and $\bar{Q}_k(t)>0$, then (\ref{eqRRdrain2})
is true as written or with the $k$ and $l$ and indices swapped. This
implies (\ref{eqRRdrain1}).

We observe that (\ref{eqIdleLimit}) is equivalent to
$
\label{eqintegralzero} \int_0^{\infty} f(\chi_m) \,d z_m= 0
$
where
\[
\chi_m:= \frac{C_i Q^{\mathbf{x}_{m}}(a_{m}t)}{a_m},\qquad z_m:=
\frac{I_i^{\mathbf{x}_{m}}(a_{m}t)}{a_m},\qquad %\label{eqfassign}
f(\cdot):= (\cdot) \wedge1.
\]
Noting that $\chi_m$ and $z_m$ meet the required conditions for
Lemma~\ref{lemintegralconv} we have,
$
\int_0^{\infty} [C_i \bar{Q}(t) ] \wedge1 \,d \bar{I}_i(t) = 0
$
which is equivalent to (\ref{eqidlelimit}).\vadjust{\goodbreak}
\end{pf*}

%%%%%%%%%%%%%%%%%%%%%%%%%%%%%%%%
%
% Theorem 3 PROOF
% 333333333333333333333333333333
%%%%%%%%%%%%%%%%%%%%%%%%%%%%%%%%%%%

%
\begin{pf*}{Proof of Theorem~\ref{thPushtoEquil}}
We first prove conclusion (i). Pick any sequence of pairs
$\{({{\mathbf{x}}_{j}},a_{j})\}$ satisfying $a_j= n_j\llVert x_j/ n_j
\rrVert_{ h \eq} \rightarrow\infty$
and $\llVert x_j/n_j \rrVert_{ h\eq}> \zeta$ for some $\zeta>0$
(a far
fluid limit sequence). To
invoke Lemma~\ref{lemsubseqtoseq}, we pick $\frak{F}$ while
simultaneously defining the process $\bar{F}(\cdot)$ according to the
expression
\[
\bar{F}(t) \triangleq\frak{F} \circ\bigl[ \bar{X}(\cdot); \bar{T}(\cdot
); \bar{
\Lambda}(\cdot); \bar{h}\bigr](t):= \bigl\llVert\bar{X}(t) \bigr
\rrVert_{ \bar{h}\eq}\qquad \forall t\geq0.
\]
Note that $\bar{F}(\llVert \bar{X}(0) \rrVert_{ \bar{h}\eq}
t)=0$ for all $t\geq t_0$
by (\ref{eqth20}), and $\frak{F}$ is easily seen to be continuous on
the topology of uniform convergence on compact sets. Since $\llVert
x_j /a_j \rrVert_{ n_jh / a_j}=1$ as argued in Corollary~\ref{cor2}, we can set
the $c$ of Lemma~\ref{lemsubseqtoseq} to~1. Applying Lemma \ref
{lemsubseqtoseq} and taking $t=t_0$ we have that
%yields
% \n{ \frac{X^\ic{l}(\n{x_l}{n_lh\eq} t) } {\n{x_l}{n_lh\eq} }}{
% for each $t\geq t_0$. Taking $t=t_0$ we have that
%
\[
\frac{1}{\llVert x_j \rrVert_{ n_jh\eq}} \bigl\llVert X^{{\mathbf
{x}}_{j }}\bigl(\llVert x_j
\rrVert_{ n_jh\eq} t_0\bigr) \bigr\rrVert_{ n_jh\eq}
\rightarrow0 \qquad\mbox{a.s.}
\]
By Lemma~\ref{lemuniformint}, $ \frac{1}{\llVert x_j \rrVert_{
n_jh\eq}}
X^{\mathbf{x}_j}(\llVert x_j \rrVert_{ n_jh\eq} t_0)$ is
uniformly integrable.
Therefore
%
%e54 #&#
%
\begin{eqnarray}
\label{eqth3seqconv} \lim_{j\rightarrow\infty} \frac{1}{\llVert x_j
\rrVert_{
n_jh\eq}}\E\bigl\llVert
X^{\mathbf{x}_j}\bigl( \llVert x_j \rrVert_{
n_jh\eq}
t_0\bigr) \bigr\rrVert_{ n_{j}h\eq} = 0.
\end{eqnarray}
%
%Applying definition (\ref{eqYdef}), and noting that our initial
%choice of sequence $\{ \ic{l} \}$ was arbitrary, up to some
%constraints, we have that the statement that follows is true. Note
%that below we have switched to scaled state notation, where recall
%$y=x/n$ and $Y^{\y}(t) = X^\ic{} (nt)/n$ for each $t$.
% \item
% For any $\const>0$, and any sequence $\{ \mathbf{x_l} = (x_l,n_l) \}$
% with $n_l \n{x_l / n_l}{h\eq} \rightarrow\infty$, and $
% %For any $\const>0$, and any sequence $\{ \mathbf{y_l} = (y_l,n_l) \}$
% %with $n_l \n{y_l}{h\eq} \rightarrow\infty$, and $
% $%\[
% %\label{eqfirstpropjfl}
% \lim_{l \rightarrow\infty} \frac{1}{\n{x_l/n_l}{h\equilib}}\E\n{
% $ %\]
% \begin{aligned}
% &\mbox{\it For any $\const>0$, and any sequence $\icb{l}$ with}
% &\mbox{\it$n_l|y_l| \rightarrow\infty$, and $|y_l| > \const$:}
% \end{aligned}\\
% \label{eqfirstpropjfl}
% \lim_{l \rightarrow\infty} \frac{1}{|y_l|}\E\left| Y^{n_l,y_l}(
%t_0|y_l|) \right|= 0.
Using that the above holds for any far fluid limit sequence, we show by
contradiction that conclusion (i) of the theorem is true.
%(\ref{eqsecondpropjfl}).
% \item
%
% For any $\const>0$, and any positive $\delta<1$ there exists
% $L_2(\const,\delta)$ such that for all $(x,n)$ satisfying $n \n{x/n}{h
%
% %For any $\const>0$, and any positive $\delta<1$ there exists
% %$L_2(\const,\delta)$ such that for all $n \n{y}{h\equilib} \geq L_2$
%and $\n{y}{h\equilib} > \const$,
% $
% \frac{1}{\n{x}{n h\equilib}}\E\n{ X^{\ic{} }( t_0 \n{x}{n h
% $
%% $ %\[
%% %\label{eqsecondpropjfla}
%% \frac{1}{\n{y}{h\equilib}}\E\n{ Y^{\y}( t_0\n{y}{h\equilib})}{h\eq}
%% $%\]
% \begin{aligned}
% &\mbox{\it For any $\const>0$, and any positive $\delta<1$ there
%exists } \nonumber\\
% &\mbox{\it$L_2(\const,\delta)$ such that for all $n|y| \geq L_2$ and
%$|y| > \const:$} \nonumber
% \end{aligned}\\
% \label{eqsecondpropjfla}
% \frac{1}{|y|}\E\left| Y^{\y}( t_0|y|) \right|\leq\delta.
Suppose conclusion~(i) %(\ref{eqsecondpropjfla})
were not true. Then for some $\zeta>0$ and some positive $\delta$,
we would
have that for any $L_2$ there would exist a pair ${{\mathbf{x}}_{}} = (x,n)$
with $\llVert x \rrVert_{ n h\equilib} \geq L_2$ and $\llVert x/n
\rrVert_{ h\equilib}>\zeta$ with
$\frac{1}{\llVert x \rrVert_{ n h\equilib}}\E\llVert X^{\mathbf
{x}}( t_0 \llVert x \rrVert_{ n h\equilib}) \rrVert_{ h\eq}
\leq\delta$.
We therefore could construct a sequence that violates (\ref
{eqth3seqconv}), which is true for any far fluid limit sequence.
A special case
of a far fluid limit sequence is when $n>L_2\zeta^{-1}$ and $\llVert
x/n \rrVert_{ h\equilib}> \zeta$. Hence we have
conclusion (i) of the theorem.

%$\y=(n,y)$
%with
%$n\n{y}{h\equilib} \geq L_2$ and $\n{y}{h\equilib}>\const$ with
%$\frac{1}{\n{y}{h\equilib}}\E\n{ Y^{\y}( t_0\n{y}{h\equilib})}{h\eq}
%> \delta$. We therefore could construct a sequence that violates
%statement (a), %(\ref{eqsecondpropjfl}),
%which is a contradiction. A special case
%of (b) %(\ref{eqsecondpropjfla})
%is when $n>L_2\const^{-1}$ and $\n{y}{h\equilib}> \const$. Hence we
%have
%conclusion (i) of the lemma.

%We begin by showing the TFL case.
We now turn to showing conclusion (ii). Pick an arbitrary sequence
of pairs $\{({{\mathbf{x}}_{j}},a_{j})\}$ satisfying $a_j={n_j}
\rightarrow
\infty$ and $\llVert x_j/n_j \rrVert_{ h\equilib} \leq\zeta$ %=
for some constant~$\zeta$ (a~near fluid limit sequence). We again invoke
Lemma~\ref{lemsubseqtoseq} by taking $\frak{F}$ to be the same
functional as before, that is,
\[
\bar{F}(t) \triangleq\frak{F} \circ\bigl[ \bar{X}(\cdot); \bar{T}(\cdot
); \bar{
\Lambda}(\cdot); \bar{h} \bigr](t):= \bigl\llVert\bar{X}(t) \bigr
\rrVert_{ \bar{h}\eq}\qquad \forall t\geq0.
\]
%
%Note that $\bar{F}(|\bar{X}(0) -\bar{h}e| t)=0 \forall t\geq
%t_0$ by Assumption (\ref{eqth30}), and $\frak{F}$ is easily seen
%to be continuous on the topology of uniform convergence on compact
%sets.
Using Lemma~\ref{lemsubseqtoseq}, and the fact that $\llVert x_j
/n_j \rrVert_{ h\equilib}\leq\zeta$ we have
$
\llVert {X^{{\mathbf{x}}_{j}}(n_jt)}/{n_j} \rrVert_{ h\eq}
\rightarrow0$
a.s.
for each $t\geq\zeta t_0$. %We are also using the fact that
%$\s{l}=n_l$ in this region of the state space.
Now take $t=t_0$,
$
\frac{1}{n_{j}}
\llVert X^{\mathbf{x}_j}(n_{j}t_0) \rrVert_{ n_{j} h \eq}
\rightarrow0$
a.s.
By Lemma~\ref{lemuniformint}, $X^{\mathbf{x}_j}(n_{j}
t_0)/n_j$ is uniformly integrable. Therefore
%
%e55 #&#
%
\begin{equation}
\label{eqth3seqconv2} \lim_{j\rightarrow\infty} \E\biggl[\frac
{1}{n_{j}} \bigl
\llVert X^{\mathbf{x}_j}(n_{j }t_0) \bigr
\rrVert_{ n_{j } h \eq} \biggr]= 0.
\end{equation}
We claim that the above implies conclusion (ii) is true by contradiction.
%Applying definition (\ref{eqYdef}), and noting that our initial
%choice of sequence $\icb{l}$ was arbitrary, up to some constraints,
%we have that the following statement is true:
% \item
% For any $\const>0$, and any sequence $\{(y_l,n_l)\}$ with $n_l
% and $\n{y_l}{h\eq} \leq\const$,
% \[
% %\label{eqfirstprop}
% \lim_{l \rightarrow\infty} \E\n{ Y^{n_l,y_l}( t_0) }{h\eq}= 0.
% \]
% \begin{aligned}
% &\mbox{\it For any $\const>0$, and any sequence $\icb{l}$ with $n_l
% &\mbox{\it and $|y_l| \leq\const$:} \nonumber
% \end{aligned}\\
% \label{eqfirstprop}
% \lim_{l \rightarrow\infty} \E\left| Y^{n_l,y_l}( t_0) \right|= 0.
%We claim that property (c) (\ref{eqfirstprop}) implies (
%We claim that fact (c) implies conclusion (ii).
% \begin{aligned}
% &\mbox{\it For any $\const>0$, and any $b>0$ there exits $L_3(
% &\mbox{\it such that for all $n\geq L_3$ and all $|y| \leq\const$:}
% \end{aligned}\\
% \label{eqsecondprop}
% \E\left| Y^{n,y}(t_0) \right| \leq b.
Suppose (ii) were not true. Then for some choice $\zeta$ and $b$, we would
have that for every constant $L_3$, there would exist an $n\geq L_3$ and
$x\dvtx  \llVert x/n \rrVert_{ h\eq} \leq\zeta$ satisfying $\E\llVert \frac
{1}{n} X^{\mathbf{x}}(n t_0) \rrVert_{ h\eq} > b$. This
would allow us to construct a sequence that violates~(\ref
{eqth3seqconv2}), which is a
contradiction.
\end{pf*}

\begin{pf*}{Proof of Lemma~\ref{lemcontract}}
%That $X$ is positive Harris recurrent follows directly from
%Theorem 3.1 of~\cite{Dai95}. The rest of the argument that follows is
%adapted from the proof of
The argument that follows is adapted from the proof of
Theorem 2.1(ii) of Meyn and Tweedie~\cite{MeynTweedieStateDep}.
We use the following fact taken from Theorem 14.2.2 of
\cite{MeynTweedie}: %to show (\ref{eqtau1lem}):
\begin{fact}[(Meyn and Tweedie~\cite{MeynTweedie})]
Suppose a discrete time Markov chain
$ \bolds{\Phi}=\{\Phi_k, k \in\mathbb{Z}^+ \}$ is defined on a
general state space $\mathrm{X}$ with transition kernel
$\Pb(x,A)=\Pb_x(\Phi_1 \in A)$, where $A\in\frak{B}(\mathrm{X})$,
the Borel subsets of $\mathrm{X}$.
If $V$ and $f$ are nonnegative measurable
functions satisfying
%{
%%
%% \label{eqMT1}
%%
%}
%
\[
\int\Pb(x,dy) V(y) \leq V(x) - f(x) + \tilde{b}1_B(x),\qquad
x\in\mathrm{X},
\]
then
\[
\E_x \Biggl[ \sum_{k=0}^{\tau_B - 1}
f( \Phi_k) \Biggr] \leq V(x) + \tilde{b},
\]
where
$
\tau_B = \inf\{k \geq1\dvtx  \Phi_k \in B\}.
$
\end{fact}
The above fact is a form of Dynkin's formula and is shown by
using the first inequality to sum bounds of the increments $E_x V(\Phi
_{k}) - E_x V(\Phi_{k+1})$ for $k\in\{0\lldots\tau_B-1\}$. Since
$1_B(\Phi_k)$ is 1 at most once for $k\in\{0,\ldots,{{\tau}_B}-1\}$ on
each sample
path, $\tilde{b}$ appears once in the final expression.

%We now turn to setting up our problem to
%make use of Fact 1. %in order to get the bound (\ref{eqtau1lem})
%on the stopping time $\tau^n_B$ as we have defined it.
We define the set $B \triangleq\{x\dvtx  \llVert x/n \rrVert_{ h\eq
}\leq\zeta\}$.
Next, we define the following functions, the first mapping each $x \in
\mathrm{X}$
to a time $m(x)$, and the second a Lyapunov function mapping each $x$
to a value:
%
%e56 #&#
%e57 #&#
%
\begin{eqnarray}
\label{eqembeddedtimes} m(x) &\triangleq& %
\cases{ n \llVert x/n
\rrVert_{ h\eq}t_0, &\quad if $x \notin B$,
\cr
%$\n{x/n}{h
n t_0, &\quad if $x \in B$,} % $\n{x/n}{h\eq} \leq\const$}
\\
% n(y) &\triangleq
% \cases{
% |y|t_0 & \mbox{ if $|y|> \zeta$}\\
% t_0 & \mbox{ if $|y|\leq\zeta$}
% } \\
\label{eqVydef}
V(x) &\triangleq&\frac{n t_0}{1-\delta}\llVert x/n \rrVert_{ h\eq}.
\end{eqnarray}
Substituting $m(x)$ for time in relation (\ref{eqlyaplem}), and
adding a term to that relation's right-hand side so that the relation
holds for
$x$ both inside and outside~$B$, we have
\begin{eqnarray*}
\E_x \biggl\llVert\frac{1}{n}X^n\bigl( m(x)
\bigr) \biggr\rrVert_{ h \eq} &\leq& \delta\llVert x/n \rrVert_{ h \eq}
+ \biggl( \sup_{ \tilde{x} \in B} \E_{\tilde{x}} \biggl\llVert
\frac
{1}{n} X^{n} ( n t_0 ) \biggr
\rrVert_{ h \eq} \biggr) 1_B(x)
\\
& \leq& \llVert x/n \rrVert_{ h \eq} + b 1_B(x)
\\
&\leq& \llVert x/n \rrVert_{ h \eq} - \frac{1-\delta}{n t_0} m(x) + (1 -
\delta+
b) 1_B(x),
\end{eqnarray*}
%
% \E_y \n{ Y^n(t_0 \n{y}{h\eq})}{h\eq}
% & \leq\delta\n{y}{h\eq} + \left(\sup_{y\in B} \E_y \n{ Y^n(t_0)}{h
%& \leq\n{y}{h\eq} - \frac{1-\delta}{t_0} n(y) + (1-\delta+
%b)1_B(y)
where the middle step follows from (\ref{eqlyaplem2}).
% \E_y|Y^{n} (n(y)) | \hspace{-30pt} \nonumber\\
% &\leq\delta|y| + \left( \sup_{y \in B} \E_y|Y^{n}( t_0)|
% &\leq|y| - \frac{1-\delta}{t_0}n(y) + %\nonumber\\
% \left( 1-\delta+ b \right)1_B(y).
By multiplying both sides by $nt_0/(1-\delta)$, we have% take $V(\cdot)
%= n t_0(1-\delta)^{-1} \n{\cdot/n }{ h \eq}$ to get
%
%e58 #&#
\begin{equation}
\label{eqhello}
\E_x\bigl[ V\bigl( X^n\bigl(m(x)
\bigr)\bigr)\bigr] \leq V(x) -m(x) + \tilde{b}1_B(x),
\end{equation}
where
%
%e59 #&#
\begin{equation}
\label{eqdefb}
\tilde{b} = nt_0 + \frac{n t_0}{1-\delta}b.
\end{equation}
%
% \label{eqhello}
% & \E_y[ V( Y^n(t_0 \n{y}{h\eq})) ] \leq V(y) -n(y) +
% &\mbox{where} \label{eqdefb} \tilde{b} = t_0 +
% \frac{t_0}{1-\delta}b.
% \label{eqhello}
% & \E_y|V\left(Y^{n}(n(y))\right)| \leq V(y) - n(y) +
% &\mbox{where} \label{eqdefb} \tilde{b} = t_0 +
% \frac{t_0}{1-\delta}b.
The transition kernel $\Pb^t$ for the Markov process $X^{n}$ is
defined by
$ %\[
\Pb^t(x,A) = \Pb_x(X^{n}(t) \in A )
$ %\]
where $A$ is any set in $\frak{B}(\mathrm{X})$, the Borel subsets
of the state space~$\mathrm{X}$. We define the discrete time
``embedded'' Markov chain $\hat{\bolds{\Phi}} = \{\hat{\Phi}_k,\break k \in
\mathbb{Z}_+\}$ with transition kernel $\hat{\Pb}$ given by
$ %\[
\hat{\Pb}(x,A) = \Pb^{m(x)}(x,A).
$ %\]
Note that
\[
\int\hat{\Pb}(x,dz)V(z) = \int\Pb^{m(x)}(x,dz) V(z) =
\E_x \bigl[ V \bigl( X^n\bigl(m(x)\bigr) \bigr) \bigr].
\]
Combining this with (\ref{eqhello}) we have
\[
\int\hat{\Pb}(x,dz)V(z) \leq V(x) - m(x) + \tilde{b} 1_B(y).
\]
Thus by fact 1,
%Recognizing this is the form of expression (\ref{eqMT1}), we may
%use Fact 1 to conclude
%
%e60 #&#
%
\begin{equation}
\label{eqnice} \E_x \Biggl[\sum_{k=0}^{\hat{\tau}_B - 1}
m(\hat{\Phi}_k) \Biggr] \leq V(x) + \tilde{b},
\end{equation}
where $\hat{\tau}_B = \inf\{k \geq1\dvtx  \hat{\Phi}_k \in B\}$.
%is the first return time of the embedded discrete time chain $
%Fix a particular $\omega$ and initial condition
%$x$.
If the embedded chain hits $B$ in $\hat{\tau}_B$ discrete
steps, then the original chain must also hit $B$ in a time less than or equal
to the sum of the embedded times. % that those discrete steps
%correspond to. It is also possible that the original chain hits
%$B$ earlier, in addition to hitting at a time equal to the sum of
%the embedded times.
Thus
\[
\inf\bigl\{t \geq0\dvtx  X^{x,n}(t) \in B \bigr\} \leq\sum
_{k=0}^{\hat{\tau}_B -
1} m(\hat{\Phi}_k) \qquad\mbox{$
\Pb_x$-a.s.}
\]
for each $x\in\mathrm{X}$. Furthermore, whenever the initial
condition $x\in B$, the first embedded time is $nt_0$ seconds by
(\ref{eqembeddedtimes}). Consequently, the time of the first
hitting of $B$ after $n t_0$ seconds expire satisfies
\[
\inf\bigl\{t \geq n t_0\dvtx  X^{x,n}(t) \in B \bigr\} \leq\sum
_{k=0}^{\hat{\tau}_B
- 1} m(\hat{\Phi}_k)
\qquad\mbox{$\Pb_x$-a.s.}
\]
for each $x\in B$. Substituting definition (\ref{eqtaulem}), taking the
expectation and using~(\ref{eqnice}), we have
\[
\E_x\bigl[\tau^n_B(n t_0)
\bigr] \leq V(x) +\tilde{b} \qquad\mbox{for all $x \in B$}.
\]
Taking the $\sup_{x\in B}$ of both sides and substituting
(\ref{eqVydef}) and (\ref{eqdefb}), we have~(\ref{eqtau1lem}).
% \sup_{y\in B} \E_y[\tau^n_B(t_0)] \leq t_0 + \frac{t_0}{1-\delta}
% \left[ \zeta+ b \right],
%which is (\ref{eqtau1lem}).
Since $B$ is closed and bounded, and arrivals are from an unbounded
distribution (\ref{equnbounded}) and spread-out (\ref{eqspreadout}),
$B$ is a petite set; see~\cite{MeynTweedie} for a discussion of petite
sets. Therefore
(\ref{eqtau1lem}) implies $X$ is positive Harris recurrent by
Theorem~4.1 of~\cite{MeynTweedie3}.
%and the fact that (\ref{equnbounded}), (\ref{eqspreadout}) make
%$B$ a petite set.
\end{pf*}
\end{appendix}

\printaddresses

\end{document}